\documentclass[12pt,leqno]{article}
\usepackage{amsmath,amsfonts}
\usepackage[dvips]{graphics,color}
\usepackage{latexsym}
\usepackage{amssymb}
\usepackage{graphicx}
\def\up{\Upsilon}
\def\upp{\Upsilon_+}
\def\upm{\Upsilon_-}
\def\lr{\log R}
\def\de{\delta}
\def\deo{\delta_{1}}
\def\k{\mathcal{K}}
\def\p{\mathcal{P}}
\def\aa{\mathcal{A}}
\def\h{\mathcal{H}}
\def\g{\mathcal{G}}
\def\n{\mathcal{N}}  
\def\l{\mathcal{L}}
\def\j{\mathcal{J}}

\def\bp{\bar{\psi}}
\def\f{\mathcal{F}}
\def\c{\mathcal{C}}
\def\sp{(s,\psi)}
\def\rsp{(\rho+s,\psi)}
\def\sbp{(1-s,\bar{\psi})}
\def\e{\mathcal{E}}
\def\d{\mathcal{D}}
\def\t{\mathcal{T}}
\def\v{\mathcal{V}}
\def\aft{A^*(\theta)}

\def\tp{\frac{1}{2\pi i}}
\def\al{\alpha}
\def\vpo{\varpi_1}
\def\vpt{\varpi_2}
\def\aht{A(\theta)}
\def\lp{\lambda_{+}}
\def\lm{\lambda_{-}}
\def\lpm{\lambda_{\pm}}
\def\u{\mathcal{U}}
\def\i{\mathcal{I}}
\def\e{\mathcal{E}}
\def\m{\mathcal{M}}
\def\q{\mathcal{Q}}
\def\r{\mathcal{R}}

\def\hh{\hat{h}}
\def\ep{\varepsilon}
\def\pt{\varphi_{\theta}}
\def\ps{\varphi_{s}}
\def\pia{\varphi_{0}}

\def\s{S(\psi)}
\def\ss{\mathcal{S}}
\def\rsp{(\rho+s,\psi)}
\def\rbp{(1-\rho-s,\bp)}
\def\rp{(\rho,\psi)}

\def\tg{\tilde{g}}

\def\tho{\vartheta\bigg(\frac{Q^{4/3}}{n}\bigg)}
\def\tht{\vartheta\bigg(\frac{Q^{2/3}}{n}\bigg)}
\def\vt{\vartheta}
\def\tr{\tilde{\mathcal{R}}}

\def\x{\mathcal{X}}
\def\y{\mathcal{Y}}

\def\w{\mathcal{W}}

\date{}
\begin{document}
\title{On the Landau-Siegel Zeros Conjecture}
\maketitle
\thanks{}
\author{\centerline {Yitang Zhang }}
\maketitle

\noindent
Table of Content \hfil
\newline

\bigskip

\noindent 
1. Introduction

2. The Set $\Psi^*$

3. Zeros of $L(s,\psi)L(s,\psi\chi)$ in $\Omega$

4. The Functions $\k_{\pm}(s,\psi)$

5. The Linear Functional $\Phi(f;\rho,\psi)$

6. Solutions to Boundary Value Problems

7. Approximation to $\pia$

8. The Fundamental Inequality: Preliminary

9. Upper Bounds for $\Theta(k)$ and $\Theta(r)$

10. Asymptotic Expression for $\Theta(\pia)$

11. The Fundamental Inequality: Completion

12. Estimation of $\e$

13. Proof of Theorem 1

References

\section *{1. Introduction}
\medskip\noindent

The main result of this paper is the following

\bigskip\noindent
 {\bf Theorem 1} {\it For any real primitive character $\chi$ of modulus $D$ we have
$$
L(1,\chi)>c_1(\log D)^{-17}(\log\log D)^{-1} \eqno (1.1)
$$
where $c_1>0$ is an effectively computable constant.}
\medskip\noindent

As a direct consequence of Theorem 1, we have

\medskip\noindent
{\bf Theorem 2} {\it Let $\chi (\mod D)$ be as in Theorem 1. Then we have

$$
L(\sigma,\chi)\ne 0\eqno (1.2)
$$
for 
$$
\sigma>1-c_2(\log D)^{-19}(\log\log D)^{-1}
$$
where $c_2>0$ is an effectively computable constant.} 
\medskip\noindent

With extra effort we can replace the factor $(\log\log D)^{-1}$ in the above results by one; even the powers of $(\log D)^{-1}$  might be slightly improved. However, the present method fails to achieve the
conjectural lower bound
$$
L(1,\chi)\gg (\log D)^{-1}.
$$
The major obstacle here is the fact that the trivial bound
$$
L'(s,\chi)\ll(\log D)^2\qquad (s\sim 1),
$$
which is applied to derive Theorem 2 from Theorem 1 and to prove Lemma 2.1, can not be essentially sharpened by present methods (see H. Iwaniec and E. Kowalski [IW, Chapter 22], for example).
\medskip\noindent

Although Theorem 2 does not completely eliminate the Landau-Siegel zeros in their original definition, our results will be sufficient for various applications in both of the analytic number theory and algebraic number theory.
\medskip\noindent

Throughout, we denote by $\chi$  a real primitive character of modulus $D$ with $D$ greater than a sufficiently large computable number. Our goal is to derive a contradiction from the following fundamental assumption
$$
L(1,\chi)<(\log D)^{-17}(\log\log D)^{-1}. \eqno (A1)
$$
\medskip\noindent

The underlying idea of our proof is initially motivated by three papers: D. Goldfeld [G], H. Iwaniec and P. Sarnak [IS] and J. B. Conrey, A. Ghosh and S. M. Gonek [CGG2]. The first paper contains asymptotic formulae (approximate functional equations) for the function $L_E(s)L_E(s,\chi)$ and its derivatives at the central point where $E$ is an elliptic curve, $L_E(s)$ is the $L$-function associated to $E$ and $L_E(s,\chi)$ is the twisted $L$-function of $E$ and $\chi$; the second paper discusses the relationship between the non-vanishing of a family of automorphic $L$-functions at the central point and the Landau-Siegel zeros; the third paper uses evaluation of discrete means over the zeros of the Riemann zeta function to derive results on the simple zeros of $\zeta(s)$. 
\medskip\noindent

Motivated by their ideas and methods, in this paper we investigate the relationship between the lower bound for $L(1,\chi)$ and the behaviors of a family of Dirichlet $L$- functions in a certain domain $\Omega$. We use $\Psi$ to denote a certain large set of Dirichlet primitive characters of which the sizes of the moduli are $\sim Q$ where
$$
Q=\exp[(\log D)^2].
$$
For technical reason we require the character $\psi\chi$ to be primitive if $\psi\in\Psi$. On the other hand, the domain $\Omega$ is chosen to be relatively small and it is not closed to the real axis. The methods used in this paper are classical, which mainly include various asymptotic formulae, the large sieve techniques and function approximations.
\medskip\noindent

The proof of Theorem 1 is based on three observations to be described below. Some of them are closely related to present researches on the Riemann zeta function, in particular, the vertical distribution of the zeros of $\zeta(s)$. Relations and applications to other problems will also be discussed. In the following descriptions we henceforth assume that (A1) holds.
\medskip\noindent

Our first observation is a "mollifier interpretation" of the assumption (A1). That is, the function 
$$
F\sp^{-1}L(s,\psi\chi)
$$
with $\psi\in\Psi$ and $s\in\Omega$, where $F\sp$ is a (relatively short) partial sum of $L\sp L(s,\psi\chi)$, behaves like a good mollifier of $L\sp$. Here the role played by the factor $F\sp^{-1}$ is, roughly speaking, to transform the real character $\chi$ into the M\"{o}bius function $\mu$. In fact, our proof is to show, by contradiction, that such a mollifier is too good to exist.
\medskip\noindent

In Section 2 and Sectino 3 we prove that the following hold for almost all of the characters $\psi$ in$\Psi$:

\medskip\noindent
(i). All the zeros of $L\sp L(s,\psi\chi)$ in $\Omega$ lie on the critical line.

\medskip\noindent
(ii). All the zeros of $L\sp L(s,\psi\chi)$ in $\Omega$ are simple.

\medskip\noindent
(iii). With a small error, the gap between any pair of consecutive zeros of $L\sp L(s,\psi\chi)$ in $\Omega$ is asymptotically equal to the average gap  $\pi(\log Q)^{-1}$.

\medskip\noindent
The proof of the above result is divided into two steps. Firstly, we formally define a subset $\Psi^*$ of $\Psi$ by three inequalities (I1), (I2) and (I3) and show that the assumption (A1) implies that almost all of the characters 
$\psi$ in $\Psi$ also belong to $\Psi^*$; secondly, we show that (i), (ii) and (iii) hold if $\psi\in\Psi^*$. It should be remarked that most of the lemmas and propositons of this paper are independent of (A1). However, without assuming (A1) they are not applicable, since, in such a situation, one is unable to determine whether the subset $\Psi^*$ is empty or not.
\medskip\noindent

Other properties of $L\sp L(s,\psi\chi)$ with $\psi\in\Psi^*$ are also obtained, among which the most remarkable one is the relation (see Lemma 3.4)
$$
F\rsp^{-1}L\rsp L(\rho+s,\psi\chi)\sim 1-Q^{-2s}, \eqno (1.3)
$$
if $\rho$ is a zero of $L\sp L(s,\psi\chi)$ in $\Omega$ and $|s|\le 3R(\log Q)^{-1}$ where
$$
R\sim\frac1{30}\log\log D.
$$

The gap assertion (iii) above is obviously related to the results and conjectures on the vertical distribution of the zeros of $\zeta(s)$. It is reasonable to believe, under the assumption of the Riemann Hypothesis, that a significant
progress toward Montgomery's pair correlation conjecture [M] or the conjectures on the small gaps between the zeros of $\zeta(s)$ (see J. B. Conrey, A. Ghosh, D. Goldston, S. M. Gonek and D. R. Heath-Brown [CGGGH], J. B. Conrey, A. Ghosh, S. M. Gonek [CGG1] and Y. Zhang [Z1]) which breaks the barrier due to the limitation of present methods, could be used here to derive a contradiction from (iii). As we are unable to do this directly, we consider the functions
$$
\k_{\pm}\sp=F\sp^{-1}L(s\pm\al,\psi) L(s\mp\al,\psi\chi)
$$
with $\psi\in\Psi^*$, where 
$$
\al=(\log Q)^{-1}.
$$
\medskip\noindent

In Section 4 we obtain some properties of the functions $\k_{\pm}\sp$, including asymptotic formulae whose proofs rely on the gap assertion (iii).
\medskip\noindent

Our second observation is a "linear functional interpretation" of the asymptotic formula for the function $\k_{+}\rsp$. For $f\in C[-1,1]$ we define
$$
\Phi(f;s,\psi)=\sum_{n\le Q^{3/2}}\frac{\lp(n)\psi(n)}{n^s}\tho f(1-\alpha\log n)
$$
$$
\qquad\qquad-\sum_{n\le Q}\frac{\lm(n)\bp(n)}{n^{1-s}}\tht f(\al\log n-1), 
$$ 
\medskip\noindent
where $\vartheta$ is a certain weight function and where $\lambda_{\pm}(n)$ denotes the $n$th coefficient in the Dirichlet series expansion of $\zeta(s\pm\al)/\zeta(s\mp\al)$. Write $\ps(x)=Q^{xs}$.
In Section 5 it is shown that if $\psi\in\Psi^*$, $\rho$ is a zero of $L\sp L(s,\psi\chi)$ satisfying $|\rho-s_0|<1$,
$$
|\theta|<3R(\log Q)^{-1},\qquad\text{and}\quad\k_+(\rho+\theta,\psi)=0,\eqno (1.4)
$$
then
$$
\Phi(\pt;\rho,\psi)\sim 0.\eqno (1.5)
$$
Consequently, it can be proved, with respect to a norm of the Sobolev type that is introduced in Section 5, that if $f\in C^1[-1,1]$ can be well approximated by a linear combination of the $\pt$ with $\theta$ satisfying (1.4), then the linear functional $\Phi(f;\rho,\psi)$ will have a small average order. This naturally leads the the problem of approximating certain functions by linear combinations of 
the $\pt$ with $\theta$ satisfying (1.4). The idea of relating asymptotic formulae for products of $L$-functions to function approximations is not only used in this paper, but also has some other applications. For example, when applied to the Riemann zeta function and combined with an asymptotic formula for the product of $\zeta(s)$ and a mollifier, it leads to some interesting results on the zeros of $\zeta(s)$ that will be given in our separate paper [Z2]. 
\medskip\noindent

Our last observation, which is the most technical part of this paper, is that if $h$ is a linear combination of the $\pt$ with $\theta$ satisfying (1.4), then one may have an asymptotic expression for $h(x)$ with $1<|x|<2$ which involves certain Dirichlet polynomials whose coefficients depend on the the function value of $h$ on $[-1,1]$ only. In  Section 7 we consider the approximation to $\pia$ on $[-1,1]$
by a linear combination of the exponential functions $\pt$ with $\theta$ satisfying (1.4) and a function $k$ which is orthogonal to the solution to a linear differential equation of Sturmian type with boundary conditions given in Section 6 (the inner product is introduced in Section 5). Such a property of $k$ is crutial in our proof. Our method actually constructs a linear combination $h$ of the $\pt$ with $\theta$ satisfying (1.4) such that, on the interval $[-1,1]$,
$$
\pia\sim h+k.
$$

Write
$$
\hh(x)=h(x+\ep)
$$
with a small positive number $\ep$ depending on $D$. Starting from Section 8, we use two different methods to calculate the quantity $\Phi(\hh;\rho,\psi)$ with $\psi\in\Psi^*$ and $\rho\in\s$ where $\s$ is a set of zeros of $L\sp L(s,\psi\chi)$ in $\Omega$. On one hand, it follows from (1.5) that
$$
\Phi(\hh;\rho,\psi)\sim 0.\eqno (1.6)
$$
On the other hand, by our third observation, we can find an asymptotic expression for $\Phi(\hh;\rho,\psi)$ which involves $\pia$. Combining this result with (1.6) we shall,
after a lengthy process of calculation in Section 8, Section 9 and Section 10 , arrive at a relation of the form
$$
R^{-1/12}\ll\e\rp \eqno (1.7)
$$
in Section 11. This formula is called "the fundamental inequality". The expression for $\e\rp$ is complicated. However, it is not difficult to estimate the discrete mean
$$
\sum_{\psi\in\Psi^*}\,\sum_{\rho\in\s}\e\rp, \eqno (1.8)
$$
by the large sieve techniques, that is done in Section 12. 
\medskip\noindent

In the last section, we prove that the expression
$$
R^{-1/12}\sum_{\psi\in\Psi^*}\,\sum_{\rho\in\s}1
$$
is of a (slightly) larger order than (1.8) if (A1) holds. This, together with (1.7), leads to a contradiction and thus completes the proof of Theorem 1.
\medskip\noindent

\section *{2. The Set $\Phi^*$}
\medskip\noindent

Write
$$\l=\log D,\qquad Q=\exp[\l^2].\eqno(2.1)
$$

\medskip\noindent
 Let $\nu(n)$ and $\upsilon(n)$ be the arithmetic functions defined by
$$
\sum_{1}^{\infty}\frac{\nu(n)}{n^s}=\zeta(s)L(s,\chi),\qquad\sum_{1}^{\infty}\frac{\upsilon(n)}{n^s}=\zeta(s)^{-1}L(s,\chi)^{-1}.\eqno(2.2)
$$

\medskip\noindent
The following lemma will play a fundamental role in the proof of Theorem 1.

\medskip\noindent
{\bf Lemma 2.1}  {\it Assume that (A1) holds. Then we have
$$
\sum_{D^5\le n\le Q^2}\frac{\nu(n)^2}{n}\ll\l^{-13}(\log\l)^{-1} \eqno (2.3)
$$
and}
$$
\sum_{D^5\le n\le D^{10}}\frac{\nu(n)^2\tau(n)^2}{n}\ll\l^{-2}(\log\l)^{-1}. \eqno (2.4)
$$

\medskip\noindent
{\it Proof}.  We prove (2.4) only, as the proof of (2.3) is analogous. Write
$$
\d(s)=\zeta(s)^{-8}L(s,\chi)^{-8}\sum_{1}^{\infty}\frac{\nu(n)^2\tau(n)^2}{n^s},
$$
so
$$
\sum_{1}^{\infty}\nu(n)^2\tau(n)^2 e^{-n/x}=\frac{1}{2\pi i}\int_{(2)}\d(s)\zeta(s)^{8}L(s,\chi)^{8}\Gamma(s)x^s\,ds.\eqno (2.5)$$
By checking the cases $\chi(p)=\pm 1,0$ respectively, we see that the factor $\d_p(s)$ in the Euler product representation
$$
\d(s)=\prod_{p}\d_p(s)
$$
satisfies
$$
\d_p(s)=
\begin{cases}
&1+O(p^{-2\sigma}),\qquad\qquad\qquad\quad\text {if}\qquad (p,D)=1,\\

&(1-p^{-4s})(1+O(p^{-2\sigma})),\qquad\text{if}\qquad p|D=1,
\end{cases}
$$
for $\sigma\ge\sigma_0>0$, the implied constant depending on $\sigma_0$ only. Hence, on the half plane $\sigma\ge 51/100$, the function $\d (s)$ is analytic and it satisfies
$$
\d (s)\ll\prod_{p|D}|1-p^{-4s}|. \eqno (2.6)
$$

Now assume $D^5\le x\le D^{10}$. Moving the line of integration to $\sigma=51/100$ and applying (2.6) and other well-known estimates, we see that, with an error $O(x^{51/100}D^{21/10})$,
the right-hand side of (2.5) is asymptotically equal to the residue of the integrand at $s=1$, that is equal to
$$
\frac{1}{2\pi i}\int_{|s-1|=\eta}\d(s)\zeta(s)^{8}L(s,\chi)^{8}\Gamma(s)x^s\,ds
$$
for any small $\eta>0$.
\medskip\noindent

We choose $\eta=\l^{-19}(\log\l)^{-1}$ and assume $|s-1|=\eta$. By (2.6) we have $\d(s)\ll 1$. On the other hand, by (A1) and the trivial bound $L'(w,\chi)\ll\l^2$ for $|w-1|\le\eta$ we get $L(s,\chi)\ll\eta\l^2$, so
$$
\zeta(s)^{8}L(s,\chi)^{8}\ll \l^{16}.
$$
(This is a crucial application of the assumption (A1)). Combining all the above estimates with (2.5) we conclude that
$$
\sum_{n\le x}\nu(n)^2\tau(n)^2\ll x\eta\l^{16}+x^{51/100}D^{21/10}\ll x\l^{-3}(\log\l)^{-1}.
$$
From this (2.4) follows by partial summation. $\Box$
\medskip\noindent

Let $\q$ denote the set of the integers $q$ ssatisfying
$$
Q<q<2Q,\qquad (6,q)=1,\qquad \mu(q)^2=1.
$$
For any $q\in\q$ let $\Psi_q$ denote the set of all the primitive characters $\psi(\mod q)$ such that $\psi\chi$ is a primitive character $(\mod [q,D])$. Let
$$
\Psi=\bigcup_{q\in\q}\Psi_q. \eqno(2.7)
$$
It will be shown in the last section that
$$
\sum_{\psi\in\Psi}1\gg Q^2.\eqno (2.8)
$$

The following result is a direct consequence of the large sieve inequality.

\medskip\noindent
{\bf Lemma 2.2} {\it For any complex numbers $a_n$ we have}
$$
\sum_{\psi\in\Psi}\bigg|\sum_{n\le Q^2} a_n\psi(n)\bigg|^2\ll Q^2\sum_{n\le Q^2}|a_n|^2.
$$
\medskip\noindent

Let $\Omega$ denote the domain
$$
|\Re(s-s_0)|<\frac{1}{2},\qquad |\Im(s-s_0)|<10\l
$$
with
$$
s_0=\frac{1}{2}+iD. \eqno(2.9)
$$
\medskip\noindent

In this and the next section we shall investigate the distribution of zeros of the function $L(s,\psi)L(s,\psi\chi)$ with $\psi\in\Psi$ in $\Omega$. Write
$$
\alpha=\l^{-2} (=(\log Q)^{-1}).\eqno(2.10)
$$
Let $\Omega_1$ and $\Omega_2$ denote the domains
$$
\Omega_1=\big\{s:\quad -2\al^{1/2}<\Re(s-s_0)<1,\qquad |\Im(s-s_0)|<1+20\l\big\},
$$

$$
\Omega_2=\big\{s:\quad -(10)^{-1}\al\log\l<\Re(s-s_0)<1/2,\qquad |\Im(s-s_0)|<1+10\l\big\}.
$$

\medskip\noindent
{\bf Lemma 2.3}  {\it  For any sequence of complex numbers $(a_n)$ write
$$
\m_1((a_n),\psi)=\sup_{s\in\Omega_1}\,\bigg\{\bigg|\sum_{n\le D^{25}}\frac{a_n\psi(n)}{n^s}\bigg|^2\bigg\},
$$
$$
\m_2((a_n),\psi)=\sup_{s\in\Omega_2}\,\bigg\{\bigg|\sum_{n\le Q^2}\frac{a_n\psi(n)}{n^s}\bigg|^2\bigg\}.
$$
Then we have

$$
\sum_{\psi\in\Psi}\m_1((a_n),\psi)\ll\l^2Q^2\sum_{n\le D^{25}}\frac{|a_n|^2}{n} \eqno (2.11)
$$
and}
$$
\sum_{\psi\in\Psi}\m_2((a_n),\psi)\ll\l^{17/5}Q^2\sum_{n\le Q^2}\frac{|a_n|^2}{n}. \eqno (2.12)
$$

\medskip\noindent
{\it Proof}  Let $\r_1$ be the rectangle with vertices at
$$
s_0-3\al^{1/2}\pm i(2+20\l),\qquad s_0+2\pm i(2+20\l).
$$

\medskip\noindent
Note that the domain $\Omega_1$ is enclosed in $\r_1$. By the Cauchy integral formula we get for any $s\in\Omega_1$
$$
\bigg|\sum_{n\le D^{25}}\frac{a_n\psi(n)}{n^s}\bigg|^2=\frac{1}{2\pi}\bigg|\int_{\r_1}\bigg(\sum_{n\le D^{25}}\frac{a_n\psi(n)}{n^w}\bigg)^2\,\frac{dw}{w-s}\bigg|.
$$ 

\medskip\noindent
Note that $n^{-2u}\ll n^{-1}$ if $n\le D^{25}$ and $u\ge 1/2-3\al^{1/2}$. Hence, by Lemma 2.2 we get for 
$w\in\r_1$
$$
\sum_{\psi\in\Psi}\bigg|\sum_{n\le D^{25}}\frac{a_n\psi(n)}{n^w}\bigg|^2\ll Q^2\sum_{n\le D^{25}}\frac{|a_n|^2}{n}.
$$

\medskip\noindent
On the other hand, note that the total length of $\mathcal{R}_1$ is $O(\l)$ and $|w-s|^{-1}<\al^{-1/2}=\l$ if $s\in\Omega_1$ and $w\in\mathcal{R}_1$. Gathering these results we obtain (2.11).
\medskip\noindent

To Prove (2.12) we choose $\r_2$ as the rectangle with vertices at
$$
s_0-\alpha\big(1+(10)^{-1}\log\l\big)\pm i(2+10\l),\qquad s_0+1\pm i(2+10\l).
$$
Note that $n^{-2u}\ll\l^{2/5}n^{-1}$ if $n\le Q^2$ and $u\ge 1/2-\al\big(1+(10)^{-1}\log\l\big)$, and $|w-s|^{-1}<\al^{-1}=\l^2$ if $s\in\Omega_2$ and 
$w\in\r_2$. The rest of the proof is analogous to that of (2.11). $\Box$
\medskip\noindent

Write
$$F(s,\psi)=\sum_{n\le D^5}\frac{\nu(n)\psi(n)}{n^s},\qquad G(s,\psi)=\sum_{n\le D^5}\frac{\upsilon(n)\psi(n)}{n^s}.
$$     
By Lemma 2.3, the number of the characters $\psi$ in $\Psi$ such that
$$
\sup_{s\in\Omega_1}\,\big\{\big|F\sp\big|\big\}\ge\frac12\l^3\log\l,
$$
 is 
 $$\ll\l^{-4}(\log\l)^{-2}Q^2\sum_{n\le D^5}\frac{\nu(n)^2}{n}\ll(\log\l)^{-2}Q^2.$$

\medskip\noindent 
 Clearly the same estimate holds with $G(s,\psi)$ in place of $F(s,\psi)$. Thus we have

\bigskip\noindent
{\bf Lemma 2.4} {\it The following inequality
$$
\sup_{s\in\Omega_1}\big\{|F(s,\psi)|+|G(s,\psi)|\big\}<\l^3\log\l\eqno (I1)
$$
holds for all but at most $O\big((\log\l)^{-2}Q^2\big)$ characters $\psi$ in $\Psi$}.
\bigskip\noindent

Now we turn to the function $F(s,\psi)G(s,\psi)-1$ which is of the form
$$
F(s,\psi)G(s,\psi)-1=\sum_{D^5<n\le D^{10}}\frac{\iota(n)\psi(n)}{n^s}
$$
with
$$|\iota(n)|\le\sum_{lm=n}|\upsilon(l)|\nu(m)\le\nu(n)\tau(n).
$$
The last inequality above, due to multiplicativity, can be checked by verifying on $n=p^l$ with $\chi(p)=\pm 1,0$ respectively. Hence, by Lemma 2.3, we see that the number of the characters $\psi$ in $\Psi$ such that
$$
\sup_{s\in\Omega_1}\{|F(s,\psi)G(s,\psi)-1|\}\ge\frac12
$$
is
$$
\ll\l^2Q^2\sum_{D^5<n\le D^{10}}\frac{\nu(n)^2\tau(n)^2}{n}.
$$

\medskip\noindent
From this and Lemma 2.1 we deduce the following

\medskip\noindent
{\bf Lemma 2.5} {\it Assume that (A1) holds. Then the following inequality
$$
\sup_{s\in\Omega_1}\big\{|F(s,\psi)G(s,\psi)-1|\big\}<\frac12 \eqno (I2)
$$
holds for all but at most $O((\log\l)^{-1}Q^2)$ characters $\psi$ in $\Psi$}.
\medskip\noindent

For any primitive character $\psi$ write
$$
\Delta(s,\psi)=\frac{L(s,\psi)}{L(1-s,\bp)}. \eqno (2.13)
$$
It is known that (see H. Davenport [D], for example)
$$
\Delta(s,\psi)=C(\psi)\bigg(\frac{k}{\pi}\bigg)^{1/2-s}\Gamma\bigg(\frac{1-s+a}{2}\bigg)\Gamma\bigg(\frac{s+a}{2}\bigg)^{-1},
$$
where $a=\big(1-\psi(-1)\big)/2$ and where $C(\psi)$ is a complex number depending on $\psi$ and satisfying $|C(\psi)|=1$.

\medskip\noindent
Suppose $\psi(\mod q)$ is in $\Psi$, so $\psi\chi$ is a primitive character $(\mod [q,D])$. Write
$$\Delta_1(s,\psi)=\Delta(s,\psi)\Delta(s,\psi\chi),$$
so
$$
L(s,\psi)L(s,\psi\chi)=\Delta_1(s,\psi)L(1-s,\bp)L(1-s,\bp\chi). \eqno (2.14)
$$

\medskip\noindent
The function $\Delta_1(s,\psi)$ ($\psi\in\Psi$) satisfies the following relations.
$$
\Delta_1(s,\psi)\Delta_1(1-s,\bp)=1,\eqno (2.15)
$$
$$
\bigg|\Delta_1\bigg(\frac12+it\bigg)\bigg|=1,\eqno (2.16)
$$
$$
\frac{\Delta_1'}{\Delta_1}\sp=-2\l^2+O(\l),\qquad (|\sigma|\le 1,\,|t-D|<21\l),\eqno (2.17)
$$
$$
|\Delta_1(s,\psi)|\le Q^{1-2\sigma}\qquad (\sigma\ge 1/2,\,t>0),\eqno (2.18)
$$
and
$$
|\Delta_1(s,\psi)|\le Q^{1-2\sigma}D^{2-4\sigma}\qquad (-2\le\sigma\le 1/2,\,|t-D|\le 20\l).\eqno (2.19)
$$

\medskip\noindent
The upper bounds (2.18) and (2.19) are crude, but they will be sufficient for our purpose.
\medskip\noindent

Let 
$$
\f(s,\psi)=F\sp+\Delta_1(s,\psi)F(1-s,\bp).\eqno (2.20)
$$

\medskip\noindent
We proceed to investigate the approximation to $L(s,\psi)L(s,\psi\chi)$ by $\f\sp$ in the domain $\Omega_2$.
To this end we introduce 
$$\varsigma(x)=\frac{1}{2\pi i}\int_{(1)}x^s\exp(s^2/4)\,\frac{ds}{s}\qquad(x>0). \eqno (2.21)
$$
Since
$$
\lim_{\epsilon\to 0^+}\varsigma(\epsilon)=0\qquad\text{and}\qquad\frac{x^s-\epsilon^s}{s}=\int_{\log\epsilon}^{\log x}e^{sy}\,dy
$$

\medskip\noindent
if $0<\epsilon\le x$, it follows, by changing the order of integration, that  
$$
\varsigma(x)=\lim_{\epsilon\to 0^+}\int_{\log\epsilon}^{\log x}\bigg\{\tp\int_{(1)}\exp(sy+s^2/4)\,ds\bigg\}\,dy\eqno (2.22)
$$
$$
=\frac{1}{\sqrt{\pi}}\int_{-\infty}^{\log x}\exp[-y^2]\,dy.
$$
Hence
$$\varsigma(x)=1+O(\exp[-(\log x)^2]),\qquad\text{if}\qquad x>1,\eqno (2.23)$$
$$
\varsigma(x)=O(\exp[-(\log x)^2]),\qquad\qquad\text{if}\qquad x\le 1.
$$
\medskip\noindent

Let $J(c)$ denote the line segment path from $c-10i\l$ to $c+10i\l$ for $c\in\mathbf{R}$. Assume that $\psi\in\Psi$ and $s\in\Omega_2$. By the residue theorem we get
$$
L\sp L(s,\psi\chi)=I_1\sp-I_2\sp+O(Q^{-21}),\eqno (2.24)
$$
where
$$
I_1\sp=\tp\int_{J(1)}L(s+w,\psi)L(s+w,\psi\chi)Q_1^w\,\exp(w^2/4)\,\frac{dw}{w},
$$
$$
I_2\sp=\tp\int_{J(-1-\sigma)}L(s+w,\psi)L(s+w,\psi\chi)Q_1^w\,\exp(w^2/4)\,\frac{dw}{w},
$$
with
$$
Q_1=Q^2D^{-10}.
$$
Here the error term arises from the integrals on the horizontal segments by routine estimation. To calculate the integral $I_1\sp$ we replace the path $J(1)$ by the line $u=1$ with a negligible error, and then apply (2.23) getting
$$
I_1\sp=F(s,\psi)+E_1\sp+O(Q^{-24}),\eqno (2.25)
$$
where
$$
E_{1}\sp=\sum_{D^5<n<Q^2}\frac{\nu(n)\psi(n)}{n^s}\varsigma\bigg(\frac{Q_1}{n}\bigg). \eqno (2.26)
$$
To calculate the integral $I_2\sp$ we apply the functional equation (2.14), the upper bound (2.10) and routine estimation to get for $w\in J(-1-\sigma)$
$$
L(s+w,\psi)L(s+w,\psi\chi)Q_1^w=\Delta_1(s+w,\psi)Q_1^w\sum_{n\le D^{25}}\frac{\nu(n)\bp(n)}{n^{1-s-w}} +O(D^{-3}).
$$
It follows that
$$
I_2\sp=I_{21}\sp+I_{22}\sp+O(D^{-3}),\eqno (2.27)
$$
where
$$
I_{21}\sp=\tp\int_{J(-1-\sigma)}\Delta_1(s+w,\psi)Q_1^w\bigg\{\sum_{n\le D^5}\frac{\nu(n)\bp (n)}{n^{1-s-w}}\bigg\}\exp(w^2/4)\,\frac{dw}{w},
$$
$$
I_{22}\sp=\tp\int_{J(-1-\sigma)}\Delta_1(s+w,\psi)Q_1^w\bigg\{\sum_{D^5<n\le D^{25}}\frac{\nu(n)\bp (n)}{n^{1-s-w}}\bigg\}\exp(w^2/4)\,\frac{dw}{w}.
$$
To calculate the integral $I_{21}\sp$ we move the path of integration from $J(-1-\sigma)$ to 
$J(2-\sigma)$ and apply the residue theorem. Then, using routine estimation of the integrals on the horizontal segments and the bound
$$
\Delta_1(s+w,\psi)Q_1^w\sum_{n\le D^5}\frac{\nu(n)\bp (n)}{n^{1-s-w}}\ll D^{-4}\qquad (w\in J(2-\sigma)),
$$
which follows from (2.18), we get
$$
I_{21}\sp=-\Delta_1(s,\psi)F(1-s,\bp)+O(D^{-4}). \eqno (2.28)
$$
To estimate the integral $I_{22}\sp$ we first note that $\Delta_1(s,\psi)\ll\l^{1/5}$ by (2.18) and (2.19). Hence, by routine estimation,
$$
I_{22}\sp=\tp\int_{J(-1-\sigma)}E_2(1-s-w,\bp)\phi(s,w)\,dw+ O(D^{-4}),\eqno (2.29)
$$
where
$$
E_2\sp=\sum_{D^5<n\le D^{25}}\frac{\nu(n)\psi (n)}{n^{s}}
$$
and
$$
\phi(s,w)=\frac{\Delta_1(s+w,\psi)Q_1^w-\Delta_1(s,\psi)}{w}\exp(w^2/4).
$$
Note that, as a function of $w$, $\phi(s,w)$ is analytic for $-1-\sigma\le u\le 1/2-\sigma$, $|v|\le 10\l$, and it satisfies
$$
\phi(s,w)\ll\frac{\l^{1/5}}{|\al+iv|}\qquad(w\in J(1/2-\sigma))
$$
by (2.18) and (2.19). On the other hand, we have
$$
|E_2(1-s-w,\bp)|=\big|E_2(1/2+i(t+v),\psi)\big|\qquad (u=1/2-\sigma).
$$
Hence, applying the Cauchy theorem to move the path of integration from $J(-1-\sigma)$ to $J(1/2-\sigma)$ and estimating the integrals on the horizontal segments trivially, we deduce from (2.29) that
$$
|I_{22}\sp|\le c_3\,\l^{1/5}\int_{-10\l}^{10\l}\big|E_2(1/2+i(t+v),\psi)\big|\,\frac{dv}{|\al+iv|}+O(D^{-4}). \eqno (2.30)
$$
where $c_3>0$ is an absolute constant.
\medskip\noindent

Now, combining (2.24), (2.25), (2.27), (2.28) and (2.30) we conclude that
$$
|L\sp L(s,\psi\chi)-\f\sp|\le |E_1\sp|+K(t,\psi)+O(D^{-3}), \eqno (2.31)
$$
for $s=\sigma+it\in\Omega_2$, where
$$
K(t,\psi)=c_3\,\l^{1/5}\int_{-10\l}^{10\l}\big|E_2(1/2+i(t+v),\psi)\big|\,\frac{dv}{|\al+iv|}.
$$
By Lemma 2.3, we see that the number of the characters $\psi$ in $\Psi$ such that 
$$
\sup_{s\in\Omega_2}\{|E_1\sp|\}\ge\frac13\l^{-24/5} 
$$
is 
$$
\ll\l^{13}Q^2\sum_{D^5\le n\le Q^2}\frac{\nu(n)^2}{n}. \eqno (2.32)
$$
On the other hand, note that if $s=\sigma+it\in\Omega_2$ and $-10\l\le v\le 10\l$, then  $1/2+i(t+v)\in\Omega_1$, so
$$
K(t,\psi)\ll\l^{1/5}\log\l\,\sup_{s\in\Omega_1}\{|E_2\sp|\}.
$$
Hence, by Lemma 2.3, the number of the characters $\psi$ in $\Psi$ such that 
$$
\sup_{|t-D|<1+20\l}\{K(t,\psi)\}\ge\frac13\l^{-24/5}
$$
is 
$$
\ll\l^{12}(\log\l)^2Q^2\sum_{D^5\le n\le D^{25}}\frac{\nu(n)^2}{n}. \eqno (2.33)
$$
By Lemma 2.1, we see that the right-hand sides of (2.32) and (2.33) are both $O((\log\l)^{-1}Q^2)$ if (A1) holds. Thus we conclude from (2.31) the following

\bigskip\noindent
{\bf Lemma 2.6} {\it Assume that (A1) holds. Then the following inequality
$$
\sup_{s\in\Omega_2}\big\{|L\sp L(s,\psi\chi)-\f\sp|\big\}<\l^{-24/5} \eqno (I3)
$$
holds for all but at most $O((\log\l)^{-1}Q^2)$ characters $\psi$ in $\Psi$}.
\medskip\noindent

Let $\Psi^*$ be the subset of $\Psi$ such that $\psi\in\Psi^*$ if and only if the inequalities (I1), (I2) and (I3) simultaneously hold. Let $\Psi^{**}$ be the complement set of $\Psi^*$ in $\Psi$. 
\medskip\noindent

From Lemma 2.4, Lemma 2.5 and Lemma 2.6 we dedrive the following

\bigskip\noindent
{\bf Proposition 2.7} {\it Assume that (A1) holds. Then we have}
$$
\sum_{\psi\in\Psi^{**}}1\ll (\log\l)^{-1}Q^2.
$$
\bigskip\noindent

We stress that the set $\Psi^*$ is formally defined; by virtue of (2.8), Proposition 2.7 shows that almost all of the characters $\psi$ in $\Psi$ belong to $\Psi^*$ under the assumption (A1) only; otherwise one is unable to determine whether $\Psi^*$ is empty or not.
\medskip\noindent

We conclude this section by proving several simple results with $\psi\in\Psi^*$ which will find application in the next two sections. 
\medskip\noindent

A direct consequence of the inequalities (I1) and (I2) is
$$
|F\sp|^{-1}<2\l^3\log\l\qquad (\psi\in\Psi^*,\,s\in\Omega_1).\eqno(2.34)
$$

\medskip\noindent
{\bf Lemma 2.8.} {\it Suppose $\psi\in\Psi^*$ and $s\in\Omega_2$. Then we have}
$$
\frac{F'}{F}\sp\ll\l\log\l.
$$

\medskip\noindent
{\it Proof.} Note that the closed disc $|w-s|\le\al^{1/2}=\l^{-1}$ is contained in $\Omega_1$, so, by (I1) and (2.34),
$$
\frac12\l^{-3}(\log\l)^{-1}<|F(w,\psi)|<\l^3\log\l,
$$
for $|w-s|\le\al^{1/2}$. Thus the lemma follows by Lemma $\al$ of E. C. Titchmarsh [T, Section 3.9]. $\Box$

\medskip\noindent
{\bf Corollary 2.9.}  {\it Suppose $\psi\in\Psi^*$, $s,\,s'\in\Omega_2$ and $|s-s'|\le(\l\log\l)^{-1}$. Then we have}
$$
\frac{F\sp}{F(s',\psi)}=1+O\big(|s-s'|\l\log\l\big).
$$

\medskip\noindent
{\it Proof.}  This follows from the relation
$$
\frac{F\sp}{F(s',\psi)}=\exp\bigg\{\int_{s'}^{s}\frac{F'}{F}(w,\psi)\,dw\bigg\}
$$
and Lemma 2.8. $\Box$
\medskip\noindent

For $\psi\in\Psi^*$ define
$$
\h\sp=\Delta_1\sp\frac{F(1-s,\bp)}{F\sp}.\eqno (2.35)
$$
By (2.17) we have
$$
|\h(1/2+it,\psi)|=1.\eqno (2.36)
$$
Furthermore, if $1/2\le\sigma<1/2+\al^{1/2}$ and $|t-D|<1+20\l$, then $1-\bar{s}\in\Omega_1$, so
$$
|F(1-s,\bp)|=|F(1-\bar{s},\psi)|\ne 0,
$$
by (2.34), and 
$$
\frac{\h'}{\h}\sp=\frac{\Delta_1'}{\Delta_1}\sp-\frac{F'}{F}\sp-\frac{F'}{F}\sbp. \eqno (2.37)
$$

\medskip\noindent
{\bf Lemma 2.10.} {\it Suppose $\psi\in\Psi^*$ and $|t-D|<1+20\l$. Then we have}
$$
\log|\h\sp|=-2\l^2\big(1+O(\l^{-1}\log\l)\big)\qquad\text{for}\quad 1/2\le\sigma<1/2+\al^{1/2},\eqno (2.38)
$$
{\it and}
$$
\h\sp\ll\l^6Q^{1-2\sigma}D^{5\sigma}\qquad\text{for}\quad \sigma\ge 1/2+\al^{1/2}.\eqno (2.39)
$$
\medskip\noindent

{\it Proof.} Assume $1/2\le\sigma<1/2+\al^{1/2}$. By (2.36) we have
$$
\log|\h\sp|=\int_{1/2}^{\sigma}\Re\bigg\{\frac{\h'}{\h}(u+it,\psi)\bigg\}\,du.
$$
Thus the formula (2.38) follows from (2.37), (2.17), Lemma 2.8 and the relation
$$
\bigg|\frac{F'}{F}\sbp\bigg|=\bigg|\frac{F'}{F}(1-\bar{s},\psi)\bigg|. 
$$
The formula (2.39) follows from (2.18), (2.34) (which trivially holds for $\sigma\ge 3/2$) and a trivial bound for $F\sbp$. $\Box$

\section  *{3. Zeros of $L(s,\psi)L(s,\psi\chi)$ in $\Omega$}
\medskip\noindent

The main result of this section is the following

\medskip\noindent
{\bf Proposition 3.1} {\it Suppose $\psi\in\Psi^*$. Then the following hold.

\medskip\noindent
(i). All the zeros of  $L(s,\psi)L(s,\psi\chi)$ in $\Omega$ lie on the critical line.

\medskip\noindent
(ii). All the zeros of $L(s,\psi)L(s,\psi\chi)$ in $\Omega$ are simple.

\medskip\noindent
(iii). If $\{\rho_1,\,\rho_2\}$ is a pair of consecutive zeros of $L(s,\psi)L(s,\psi\chi)$ in $\Omega$,
then}
$$
|\rho_1-\rho_2|=\pi\al+O\big(\l^{-3}(\log\l)^2\big).
$$

\medskip\noindent

The proof of Proposition 3.1 is obtained from Lemma 3.2, Lemma 3.3 and Lemma 3.5 below.
\medskip\noindent

We henceforth assume that $\psi\in\Psi^*$ in this section. Instead of dealing with the function $L(s,\psi)L(s,\psi\chi)$ directly, we consider the function
$$
\g\sp=\frac{L(s,\psi)L(s,\psi\chi)}{F\sp} \eqno (3.1)
$$
which has the same zeros as $L(s,\psi)L(s,\psi\chi)$ in $\Omega_1$ by (2.34). The advantage is that the function $\g\sp$ is directly related to $\h\sp$ by (I3).
\medskip\noindent

Note that $\Omega_2\subset\Omega_1$ and
$$
\frac{\f\sp}{F\sp}=1+\h\sp.
$$
It follows from (I3) and (2.34) that
$$
\g\sp=1+\h\sp+O(\l^{-9/5}\log\l)\qquad (s\in\Omega_2). \eqno (3.2)
$$
Note that if $1/2-\al<\sigma<1$ and $|t-D|<20\l$, then the closed disc $|w-s|\le\al$ is contained in $\Omega_2$. Hence, by (3.2) and the Cauchy integral formula we get 
$$
\g'\sp=\h'\sp+O(\l^{1/5}\log\l)\qquad(1/2-\al<\sigma<1,\quad |t-D|<20\l). \eqno (3.3)
$$

\medskip\noindent
From these relations we can derive the assertions (i) and (ii) of Proposition 3.1.

\medskip\noindent
{\bf Lemma 3.2} {\it Suppose $\psi\in\Psi^*$, $1/2<\sigma<1$ and $|t-D|<20\l$. Then we have}
$$
\g\sp\neq 0.
$$

\medskip\noindent
{\it Proof.} In the case $1/2+\al\le\sigma<1$ we have $|\h\sp|\le e^{-1}$ by Lemma 2.10, so the assertion follows from 
(3.2). In the case $1/2<\sigma<1/2+\al$ we have, by (3.3),
$$
\g(1-\bar{s},\psi)-\g\sp=\h(1-\bar{s},\psi)-\h\sp+O\big((2\sigma-1)\l^{1/5}\log\l\big).\eqno (3.4)
$$
On the other hand, by (2.15) we get
$$
|\h\sp|^{-1}=|\h(1-\bar{s},\psi)|.\eqno (3.5)
$$
Hence, by Lemma 2.10,
$$
|\h(1-\bar{s},\psi)-\h\sp|\gg(2\sigma-1)\l^2. \eqno (3.6)
$$
From (3.4) and (3.6) it follows that $\g(1-\bar{s},\psi)-\g\sp\neq 0$, which implies that $\g\sp\neq 0$ by the symmetry of the non-trivial zeros of $L(s,\psi)L(s,\psi\chi)$ with respect to the critical line. $\Box$

\medskip\noindent
{\bf Lemma 3.3.} {\it Suppose $\psi\in\Psi^*$ and $|t-D|<10\l$. Then we have}
$$
|\g'(1/2+it,\psi)|\gg\l^2.
$$

\medskip\noindent
{\it Proof.}  By (2.37), (2.17) and Lemma 2.8 we get
$$
\bigg|\frac{\h'}{\h}(1/2+it,\psi)\bigg|\gg\l^2.
$$
Combining this with (2.36) and (3.3) we get the assertion. $\Box$
\medskip\noindent

Lemma 3.2 and Lemma 3.3 together imply the assertions (i) and (ii) of Proposition 3.1. To prove the gap assertion (iii) of Proposition 3.1 we need a result which will also play an important role in the rest of this paper.

\medskip\noindent
{\bf Lemma 3.4.} {\it Suppose $\psi\in\Psi^*$, $\rho$ is a zero of $L(s,\psi)L(s,\psi\chi)$ in $\Omega$, and $|s|<(10)^{-1}\al\log\l$. Then we have}
$$
\g\rsp=1-Q^{-2s}+O\big(\l^{-1}(\log\l)^2Q^{-2\sigma}\big).
$$

\medskip\noindent
{\it Proof.} It follows from (3.2) that 
$$
\h\rp=-1+O(\l^{-9/5}\log\l).\eqno (3.7)
$$
On the other hand, by (2.38), (2.18) and Lemma 2.8 we get
$$
\frac{\h\rsp}{\h\rp}=\exp\bigg\{\int_{0}^{s}\frac{\h'}{\h}(\rho+w,\psi)\,dw\bigg\}\eqno (3.8)
$$
$$
$$
$$
\qquad\qquad\qquad\qquad=\exp\big\{-2\l^2s+O(\l^{-1}(\log\l)^2)\big\}
$$
$$
$$
$$
\qquad\qquad\qquad=Q^{-2s}\big\{1+O(\l^{-1}(\log\l)^2)\big\}.
$$

\medskip\noindent
Note that $\l^{-9/5}\log\l<\l^{-1}(\log\l)^2Q^{-2\sigma}$. From (3.7) and (3.8) it follows that
$$
\h\rsp=-Q^{-2s}+O\big\{\l^{-1}(\log\l)^2Q^{-2\sigma}\big\}.\eqno (3.9)
$$
By (3.9) and (3.2) we complete the proof. $\Box$
\medskip\noindent

The gap assertion (iii) of Proposition 3.1 is deduced from the following

\medskip\noindent
{\bf Lemma 3.5.} {\it Suppose $\psi\in\Psi^*$ and $\rho$ is a zero of $L(s,\psi)L(s,\psi\chi)$ in $\Omega$. Then the following hold.

\medskip\noindent
(a). $\g(\rho+it,\psi)\neq 0$ for $0<t\le\pi\al-c_4\l^{-3}(\log\l)^2$,
where $c_4>0$ is a sufficiently large constant.

\medskip\noindent
(b). The function $\g\rsp$ has a zero in the disc} $|s-\pi i\al|< c_4\l^{-3}(\log\l)^2$.

\bigskip\noindent
{\it Proof.}  Assume $0<t\le\pi\al-c_4\l^{-3}(\log\l)^2$. Then, by (3.3),
$$
\g(\rho+it,\psi)=\int_{0}^{t}\h'(\rho+iv,\psi)\,dv+O(\l^{1/5}(\log\l)t).
$$
For $0\le v\le t$ we have
$$
\h(\rho+iv,\psi)=-Q^{-2iv}+O\big(\l^{-1}(\log\l)^2\big)
$$
by Lemma 3.4 and (3.2), and
$$
\frac{\h'}{\h}(\rho+iv,\psi)=-2\l^2+O(\l\log\l)
$$
by (2.37), (2.17) and Lemma 2.8. It follows that
$$
\g(\rho+it,\psi)=2\l^2\int_{0}^{t}Q^{-2iv}\,dv+O\big(\l(\log\l)^2t\big)=i(Q^{-2it}-1)+O\big(\l(\log\l)^2t\big).\eqno (3.10)
$$
Recall that $\al=\l^{-2}$. Since
$$
|Q^{-2it}-1|=2\sin(\l^2t)\ge
\begin{cases} (4/\pi)\l^2t,\qquad&\text{if}\quad 0<t\le\pi\al/2,\\
\,\\
(4c_4/\pi)\l^{-1}(\log\l)^2,\qquad&\text{if}\quad \pi\al/2<t\le\pi\al-c_4\l^{-3}(\log\l)^2,
\end{cases}
$$
the right-hand side of (3.10) can not be equal to zero as the constant $c_4$ is sufficiently large. This proves (a).
\medskip\noindent

To prove (b) we note that $1-Q^{-2\pi i\al}=0$ and
$$
|1-Q^{-2s}|=2c_4\l^{-1}(\log\l)^2+O\big(\l^{-2}(\log\l)^4\big),
$$
if $|s-\pi i\al|=c_4\l^{-3}(\log\l)^2$, the implied constant depending on $c_4$. It follows, by Lemma 3.4, that
$$
|\g\rsp|>|\g(\rho+\pi i\al,\psi)|\qquad\text{for}\quad |s-\pi i\al|=c_4\l^{-3}(\log\l)^2.
$$
Thus, by the maximum modulus principle we get (b), and complete the proof of Proposition 3.1. $\Box$
\medskip\noindent

We conclude this section by proving some consequences of Proposition 3.1 which will find application in the next section.

\medskip\noindent
{\bf Lemma 3.6.} {\it Suppose $\psi\in\Psi^*$, $0<\sigma<1$ and $|t-D|<20\l-1$ such that
$$
L(s+w,\psi)L(s+w,\psi\chi)\ne 0,\eqno (3.11)
$$
 for any $w=u+iv$ with 
$$
-(1/10)\al<u<(21/20)\al\qquad\text{and}\quad |v|<(1/10)\al.
$$
 Then we have}
$$
1\ll\frac{|L(s+2\al,\psi)|}{|L\sp|}\ll 1.
$$

\medskip\noindent
{\it Proof.} We have
$$
\log\frac{|L(s+2\al,\psi)|}{|L\sp|}=\int_{0}^{2\al}\Re\bigg\{\frac{L'}{L}(s+u,\psi)\bigg\}\,du.
$$

\medskip\noindent
 If $\rho$ is a zero of $L\sp$ such that $|\Im\{\rho\}-t|<1$, then $\rho\in\Omega$, so $\Re\{\rho\}=1/2$ by Proposition 3.1. Write $\rho=1/2+i\gamma$. We have
$$
\Re\bigg\{\frac{L'}{L}(s+u,\psi)\bigg\}=\sum_{\substack{\rho\\|\gamma-t|<1}}\frac{\sigma+u-1/2}{(\sigma+u-1/2)^2+(\gamma-t)^2}+O(\l^2)\qquad(0\le u\le 2\al),
$$
where $\rho$ runs through the zeros of $L\sp$. Note that $|s+u-\rho|\gg\al$ for any zero $\rho$ of $L\sp$ by (3.11). Hence, by Proposition 3.1 (iii), the above sum is also $O(\l^2)$. This completes the proof. $\Box$

\medskip\noindent
{\bf Lemma 3.7.}  {\it Suppose $\psi\in\Psi^*$, $\rho=1/2+i\gamma$ is a zero of $L\sp$ in $\Omega$ such that $|\gamma-D|<20\l-2$. Then we have}
$$
\frac{L(\rho\pm 2\al,\psi)}{L'\rp}\ll\al. \eqno (3.12)
$$

\medskip\noindent
{\it Proof.} By Proposition 3.1 and the residue theorem, the left-hand side of (3.12) is equal to
$$
\tp\int_{|s|=3\al}\frac{L(\rho+s\pm 2\al,\psi)}{L\rsp}\,ds.
$$
Moreover, by Proposition 3.1 and Lemma 3.6 we get 
$$
\frac{L(\rho+s\pm 2\al,\psi)}{L\rsp}\ll 1
$$
 for $|s|=3\al$. This completes the proof. $\Box$

\section *{4. The Functions $\k_{\pm}\rsp$}
\medskip\noindent

For $\psi\in\Psi^*$ let $\s$ denote the set of zeros of $L(s,\psi)L(s,\psi\chi)$ in the strip $|t-D|<1$. In this and the next seven sections (except for Section 6) we shall assume that 
$$
\psi\in\Psi^*\qquad\text{and}\qquad\rho\in\s,\eqno (4.1)
$$
with the aim of deriving the fundamental inequality described in Section 1. At many places in these sections we shall not repeat the assumption (4.1) that applies throughout.
\medskip\noindent

Write
$$
R=\pi\bigg[\frac{\log\l}{30\pi}\bigg]-\frac{\pi}{2},\qquad\omega=\al R.\eqno (4.2)
$$

\medskip\noindent
In this section, we investigate the behaviors of the functions $\k_{\pm}\rsp$ with $|s|\le 3\omega$, where
$$
\k_+\sp=\frac{L(s+\al,\psi) L(s-\al,\psi\chi)}{F\sp},\qquad\k_-\sp=\frac{L(s-\al,\psi) L(s+\al,\psi\chi)}{F\sp}.\eqno (4.3)
$$
Recall that the function $\g\sp$ is defined by (3.1). We have
$$
\k_+\rsp=\frac{F(\rho+s-\al,\psi)}{F\rsp}\frac{L(\rho+s+\al,\psi)}{L(\rho+s-\al,\psi)}\g(\rho+s-\al,\psi),\eqno (4.4)
$$
and
$$
\k_-\rsp=\frac{F(\rho+s+\al,\psi)}{F\rsp}\frac{L(\rho+s-\al,\psi)}{L(\rho+s+\al,\psi)}\g(\rho+s+\al,\psi).\eqno (4.5)
$$
Recall that the function $\Delta\sp$ is defined by (2.13). We have
$$
\k_+\rsp=\frac{F\rbp}{F\rsp}\,\Delta(\rho+s+\al,\psi)\Delta(\rho+s-\al,\psi\chi)\k_-\rbp,\eqno (4.6)
$$
and
$$
\k_-\rsp=\frac{F\rbp}{F\rsp}\,\Delta(\rho+s-\al,\psi)\Delta(\rho+s+\al,\psi\chi)\k_+\rbp.\eqno (4.7)
$$

\medskip\noindent
Note that $Q^{\al}=e$ and $Q^{-2\sigma}\le e^{2R}<\l^{1/5}$ if $|s|\le\omega$. From (4.4), (4.5), Lemma 3.4 and Corollary 2.9 we derive the following

\medskip\noindent
{\bf Lemma 4.1.} {\it Suppose $|s|\le 3\omega$. Then we have}
$$
\k_+\rsp\k_-\rsp=(1-e^2Q^{-2s})(1-e^{-2}Q^{-2s})+O(\l^{-3/5}R^2).
$$
\medskip\noindent
Here we have written $R$ for the factor $\log\l$ in the error terms in Lemma 3.4 and Corollary 2.9.

By a familiar property of $\Delta\sp$ we have
$$
\Delta(\rho+s\pm\al,\psi)\Delta(\rho+s\mp\al,\psi\chi)\ll Q^{-2\sigma}\qquad (|s|\le 3\omega).\eqno (4.8)
$$

\medskip\noindent
{\bf Lemma 4.2.} {\it Suppose $|s|\le 3\omega$. Then we have
$$
\k_{\pm}\rsp\ll 1 \eqno(4.9)
$$
if $\sigma\ge 0$, and 
$$
Q^{2s}\k_{\pm}\rsp\ll 1 \eqno (4.10)
$$
if} $\sigma\le 0$.

\medskip\noindent
{\it Proof.} First we prove (4.9). By the maximum modulus principle, it suffices to show that (4.9) holds for
$$
\sigma\ge -2\al,\,|s|=\omega,\qquad\text{or}\qquad\sigma=-2\al,\,|t|\le\omega.
$$
By Lemma 3.6 we see that 
$$
1\ll\bigg|\frac{L(\rho+s+\al,\psi)}{L(\rho+s-\al,\psi)}\bigg|\ll 1
$$
if $\sigma=-2\al$ and $|t|\le\omega$ or $|s|=\omega$.
On the other hand, by Lemma 3.4 we have
$$
\g(\rho+s\pm\al,\psi)\ll 1\qquad\text{for}\qquad\sigma\ge -2\al,\,|s|\le\omega.
$$
Combining these estimates with (4.4), (4.5) and Corollary 2.9 we get (4.9).
\medskip\noindent

Next we prove (4.10). Assume $\sigma\le 0$ and $|s|\le\omega$. Note that
$$
|\k_{\pm}\rbp|=|\k_{\pm}(\rho-\bar{s},\psi)|.
$$
Hence, by (4.6), (4.7), (4.8) and Corollary 2.9  we get
$$
\k_{\pm}\rsp\ll Q^{-2\sigma}.
$$
This proves (4.10). $\Box$
\medskip\noindent

Let $\c$ denote the circle $|s|=\omega$.  Let 
$$
T\rp=T_1\rp\cup T_2\rp,
$$ 
where $T_1\rp$ denotes the set of all the zeros of $\k_+\rsp$ inside $\c$ and where $T_2\rp$ denotes the set of all the zeros of $L(\rho+s+\al,\psi)$ in the domain $\omega<|s|<3\omega$. By Proposition 3.1, each $\theta\in T\rp$ is of the form
$$
\theta=\pm\al+\pi il\al+O(\l^{-3}R^3),\eqno (4.11)
$$
with
$$
l\in\mathbf{Z}\qquad\text{and}\qquad|l|\le\frac{\log\l}{10\pi};
$$
conversely, for any $l$ satisfying the above relations, there exists at most one $\theta\in T\rp$ such that either
$$
\theta=\al+\pi il\al+O(\l^{-3}R^3)\qquad\text{or}\qquad\theta=-\al+\pi il\al+O(\l^{-3}R^3).\eqno (4.12)
$$
Furthermore, by Proposition 3.1 we have
$$
1-e^{-2}Q^{-2\theta}\ll\l^{-1}R^3\qquad\text{if}\quad\theta\in T_2\rp.\eqno (4.13)
$$

\medskip\noindent

Our last task of this section is to establish asymptotic formulae for $\k_{\pm}\rsp$ that essentially rely on the gap assertion (iii) of Proposition 3.1. We introduce
$$
\vt(x)=\tp\int_{(1)}x^sY(s)\,ds\qquad(x>0), \eqno (4.14)
$$
where
$$
Y(s)=\frac{5\al(Q^{s/10}-Q^{-s/10})\exp(s^2/4)}{s^2}.
$$
Since
$$
\frac{\al(Q^{s/10}-Q^{-s/10})}{s}=\int_{-1/10}^{1/10}Q^{ys}\,dy,
$$
it follows that
$$
\vt(x)=5\int_{-1/10}^{1/10}\varsigma(Q^yx)\,dy 
$$

\medskip\noindent
with the function $\varsigma(x)$ given by (2.21). Hence, by (2.22) and (2.23) we get
$$
0<\vt(x)<1, \eqno (4.15)
$$
\medskip\noindent
$$
\vt(x)=1+\exp[-(\log x-\l^2/10)^2]\quad\qquad\text{if}\quad x>Q^{1/10},\eqno (4.16)
$$
$$
\quad\vt(x)\ll\exp[-(\log x+\l^2/10)^2]\qquad\qquad\text{if}\quad x<Q^{-1/10}.
$$
\medskip\noindent

Let $\lp(n)$ and $\lm(n)$ be the arithmetic functions defined by
$$
\sum_{1}^{\infty}\frac{\lp(n)}{n^s}=\frac{\zeta(s+\al)}{\zeta(s-\al)},\qquad
\sum_{1}^{\infty}\frac{\lm(n)}{n^s}=\frac{\zeta(s-\al)}{\zeta(s+\al)}.\eqno (4.17)
$$

\medskip\noindent
{\bf Lemma 4.4} {\it Suppose $|s|\le 3\omega$. Then we have
$$\k_+\rsp=\sum_{n\le Q^{3/2}}\frac{\lp(n)\psi(n)}{n^{\rho+s}}\tho\eqno (4.18)
$$
$$
\qquad\qquad\qquad\qquad\qquad\qquad-Q^{-2s}\sum_{n\le Q}
\frac{\lm(n)\bp(n)}{n^{1-\rho-s}}\tht+O(\l^{-3/5}),
$$
and}
$$
\k_-\rsp=\sum_{n\le Q^{3/2}}\frac{\lm(n)\psi(n)}{n^{\rho+s}}\tho\eqno (4.19)
$$
$$
\qquad\qquad\qquad\qquad\qquad\qquad-Q^{-2s}\sum_{n\le Q}
\frac{\lp(n)\bp(n)}{n^{1-\rho-s}}\tht+O(\l^{-3/5}).
$$

\medskip\noindent
{\it Proof.} We prove (4.18) only, as the proof of (4.19) is analogous. By the maximum modulus principle, it suffices
to prove (4.19) with $|s|=3\omega$ which is thus assumed, so
$$
L(\rho+s-\al,\psi)\ne 0,\eqno (4.20)
$$
and, by Lemma 3.6,
$$
\frac{L(\rho+s+\al,\psi)}{L(\rho+s-\al,\psi)}\ll 1.\eqno (4.21)
$$

\medskip\noindent
Write
$$
\sigma^*=\max\{2\al,\,\sigma+2\al\}.
$$

\medskip\noindent
Let $Z$ denote the set of all the zeros of $L(w,\psi)$ in the strip $|w-(\rho+s)|<10\l-3$. Note that $Z$ depends on $\rho+s$ and $\psi$ and $Z\subset\Omega$. Also, note that the function $Y(w)$ has only a simple pole at $w=0$ and the residue is equal to one. Write $w=u+iv$. Since 
$$
\Re\{\rho+s+w-\al\}<1/2\qquad\text{if}\quad u=-\sigma^*, 
$$
by Proposition 3.1, the residue theorem and (4.21), there exists a rectangle $\r^*$ with vertices at
$$
-\sigma^*\pm i\big(10\l+O(1)\big),\qquad 1\pm i\big(10\l+O(1)\big)
$$
such that
$$
\frac{L(\rho+s+w+\al,\psi)}{L(\rho+s+w-\al,\psi)}\ll 1\eqno (4.22)
$$

\medskip\noindent
if $w$ lies on the horizontal segments of $\r^*$ (see Lemma 3.6), and such that
$$
\tp\int_{\r^*}\frac{L(\rho+s+w+\al,\psi)}{L(\rho+s+w-\al,\psi)}\kappa(s,w)Y(w)\,dw=\frac{L(\rho+s+\al,\psi)}{L(\rho+s-\al,\psi)}\kappa(s,0)\eqno (4.23)
$$
$$
\qquad+\sum_{\rho^*\in Z}\frac{L(\rho^*+2\al,\psi)}{L'(\rho^*,\psi)}\kappa(s,\rho^*-\rho-s+\al)Y(\rho^*-\rho-s+\al)
$$
where
$$
\kappa(s,w)=Q^{(4/3)w}-e^2Q^{-2s}Q^{-(2/3)w}. 
$$
Recall that $J(c)$ is the line segment path from $c-10i\l$ to $c+10i\l$. We claim that
$$
\k_+\rsp=\n_+-\n_-+O(\l^{-3/5}),\eqno (4.24)
$$
where
$$
\n_+=\tp\int_{J(1)}\frac{L(\rho+s+w+\al,\psi)}{L(\rho+s+w-\al,\psi)}\kappa(s,w)Y(w)\,dw,
$$
$$
$$
$$
\n_-=\tp\int_{J(-\sigma^*)}\frac{L(\rho+s+w+\al,\psi)}{L(\rho+s+w-\al,\psi)}\kappa(s,w)Y(w)\,dw.
$$
Since 
$$
\kappa(s,w)Y(w)\ll Q^{-23}
$$
if $-\sigma^*\le u\le 1$ and $|v|=10\l+O(1)$, it follows from (4.22) that
$$
\tp\int_{\r^*}\frac{L(\rho+s+w+\al,\psi)}{L(\rho+s+w-\al,\psi)}\kappa(s,w)Y(w)\,dw=\n_+-\n_-+O(Q^{-23}),\eqno (4.25)
$$
On the other hand, by Lemma 3.4 we get
$$
\kappa(s,0)=\g(\rho+s-\al,\psi)+O\big(\l^{-4/5}R^2\big). 
$$
Combining this result with (4.21), (4.4) and Corollary 2.9 we get
$$
\frac{L(\rho+s+\al,\psi)}{L(\rho+s-\al,\psi)}\kappa(s,0)=\k_+\rsp+O\big(\l^{-4/5}R^2\big). \eqno (4.26)
$$
Now, by (4.23), (4.25), (4.26) and Lemma 3.7, to complete the proof of (4.24), it suffices to prove that
$$
\sum_{\rho^*\in Z}\big|\kappa(s,\rho^*-\rho-s+\al)Y(\rho^*-\rho-s+\al)\big|\ll\l^{7/5}. \eqno (4.27)
$$

\medskip\noindent
Assume $\rho^*\in Z$. By direct calculation we get 
$$
\kappa(s,\rho^*-\rho-s+\al)=Q^{(4/3)(\rho^*-\rho-s+\al)}\big(1-Q^{-2(\rho^*-\rho)}\big),
$$
so
$$
\kappa(s,\rho^*-\rho-s+\al)Y(\rho^*-\rho-s+\al)\ll\frac{\al\big(1+Q^{-(43/30)\sigma}\big)\big|1-Q^{-2(\rho^*-\rho)}\big|}{\big|\rho^*-\rho-s+\al\big|^2}.\eqno (4.28)
$$

\medskip\noindent
Write $\al^*=\pi\l^{-6/5}$. By Proposition 3.1, if $\rho^*\in Z$ and $|\rho^*-\rho|<\al^*$, then there exists an integer $l$ such that $|l|\le\l^{4/5}$ and
$$
\rho^*-\rho=\pi il\al+O(|l|\l^{-3}R^2),\eqno (4.29)
$$
so
$$
1-Q^{-2(\rho^*-\rho)}\ll |l|\l^{-1}R^2;
$$
conversely, if $l\in\mathbf{Z}$ and $|l|\le\l^{4/5}$, then there exists at most one $\rho^*\in Z$ satisfying (4.29).
Note that $|\pi il\al-s+\al|\gg\al$ for any $l\in\mathbf{Z}$ since $|s|=3\omega$. It follows that
$$
\sum_{\substack{\rho^*\in Z\\|\rho^*-\rho|<\al^*}}\frac{\big|1-Q^{-2(\rho^*-\rho)}\big|}{\big|\rho^*-\rho-s+\al\big|^2}\ll \l^{-1}R^2\sum_{|k|\le\l^{4/5}}\frac{|k|}{|\pi ik\al-s+\al|^2}\ll\l^3R^3.\eqno (4.30)
$$

\medskip\noindent
On the other hand, by Proposition 3.1 and routine estimation we get
$$
\sum_{\substack{\rho^*\in Z\\|\rho^*-\rho|\ge\al^*}}\frac{\big|1-Q^{-2(\rho^*-\rho)}\big|}{\big|\rho^*-\rho-s+\al\big|^2}\ll\l^4\sum_{k\ge\l^{4/5}}\frac{1}{k^2}\ll\l^{16/5}.\eqno (4.31)
$$
Note that 
$$
Q^{-(43/30)\sigma}<\l^{43/300}<\l^{1/5}.
$$
By (4.28), (4.30) and (4.31) we get (4.27), and complete the proof of (4.24).
\medskip\noindent

It now remains to calculate the integrals $\n_+$ and $\n_-$. Note that
$$
\frac{L(\rho+s+w+\al,\psi)}{L(\rho+s+w-\al,\psi)}=\sum_{1}^{\infty}\frac{\lp(n)\psi(n)}{n^{\rho+s+w}}\qquad (u=1).
$$
To calculate $\n_+$ we can replace the path $J(1)$ by the line $u=1$ and replace the factor $\kappa(s,w)$ by $Q^{(4/3)w}$ with negligible errors. It follows that
$$
\n_+=\sum_{1}^{\infty}\frac{\lp(n)\psi(n)}{n^{\rho+s}}\tho+O(Q^{-1/2}).\eqno (4.32)
$$
To calculate $\n_-$ we first appeal to (2.13), Lemma 3.6 and the relation
$$
\frac{\Delta(\rho+s+w+\al,\psi)}{\Delta(\rho+s+w-\al,\psi)}=Q^{-2\al}+O(\l^{-1}),\qquad\big(w\in J(-\sigma^*)\big)
$$
to get
$$
\frac{L(\rho+s+w+\al,\psi)}{L(\rho+s+w-\al,\psi)}=e^{-2}\frac{L(1-\rho-s-w-\al,\bp)}{L(1-\rho-s-w+\al,\bp)}+O(\l^{-1}),\qquad\big(w\in J(-\sigma^*)\big),
$$
the contribution from the error term to $\n_-$ being $O(\l^{-4/5})$ by routine estimation. Hence
$$
\n_-=\frac{e^{-2}}{2\pi i}\int_{J(-\sigma^*)}\frac{L(1-\rho-s-w-\al,\bp)}{L(1-\rho-s-w+\al,\bp)}\kappa(s,w)Y(w)\,dw+O(\l^{-4/5}).
$$

\medskip\noindent
Next we apply the Cauchy theorem to move the path of integration from $J(-\sigma^*)$ to $J(-1)$ with negligible errors arising from the integrals on the horizontal segments. Then, with negligible errors, we can replace the path $J(-1)$ by the line $u=-1$ and replace the factor 
$e^{-2}\kappa(s,w)$
by $-Q^{-2s}Q^{-(2/3)w}$ to get

\medskip\noindent

$$
\n_-=-\frac{Q^{-2s}}{2\pi i}\int_{(-1)}\frac{L(1-\rho-s-w-\al,\bp)}{L(1-\rho-s-w+\al,\bp)}Q^{-(2/3)w}Y(w)\,dw+O(\l^{-4/5}).\eqno (4.33)
$$
Since
$$
\frac{L(1-\rho-s-w-\al,\bp)}{L(1-\rho-s-w+\al,\bp)}=\sum_{1}^{\infty}\frac{\lm(n)\bp(n)}{n^{1-\rho-s-w}}\qquad(u=-1),
$$
and, by a change of variable,
$$\tp\int_{(c)}x^wY(w)\,dw=-\vt(1/x)\qquad\text{for}\quad c<0,
$$
it follows from (4.33) that
$$
\n_-=Q^{-2s}\sum_{1}^{\infty}\frac{\lm(n)\bp(n)}{n^{1-\rho-s}}\tht+O(\l^{-4/5}). \eqno (4.34)
$$
Finally, by virtue of (4.16), we see that the terms in the sum in (4.32) with $n>Q^{3/2}$ and the terms in the sum in (4.34) with $n>Q$ can be neglected. Hence, combining (4.32) and (4.34) with (4.24) we get (4.18). $\Box$

\section  *{5. The Linear Functional $\Phi(f;\rho,\psi)$}
\medskip\noindent

 For $f\in C[-1,1]$ we define the linear functional
$$
\Phi(f;\rho,\psi)=\sum_{n\le Q^{3/2}}\frac{\lp(n)\psi(n)}{n^{\rho}}\tho f(1-\alpha\log n)\eqno (5.1)
$$
$$
\qquad\qquad-\sum_{n\le Q}\frac{\lm(n)\bp(n)}{n^{1-\rho}}\tht f(\al\log n-1). 
$$ 
Let
$$\ps(x)=Q^{xs}.$$
Note that
$$\ps(1-\alpha\log n)=\bigg(\frac{Q}{n}\bigg)^s,\qquad \ps(\alpha\log n-1)=\bigg(\frac{Q}{n}\bigg)^{-s},
$$
so
$$\Phi(\ps;\rho,\psi)=Q^s\bigg\{\sum_{n\le Q^{3/2}}\frac{\lp(n)\psi(n)}{n^{\rho+s}}\tho-Q^{-2s}\sum_{n\le Q}\frac{\lm(n)\psi(n)}{n^{1-\rho-s}}\tht\bigg\}.$$

\medskip\noindent
Hence, by Lemma 4.3 and the inequality $|Q^s|<\l^{1/10}$ for $|s|\le 3\omega$ we derive the following 

\medskip\noindent
{\bf Lemma 5.1} {\it Suppose $|s|\le 3\omega$. Then we have
$$
\Phi(\ps;\rho,\psi)=Q^s\k_+\rsp+O(\l^{-1/2}).
$$

\medskip\noindent
In particular, if $\theta\in T\rp$, then}
$$
\Phi(\pt;\rho,\psi)\ll\l^{-1/2}.
$$
\medskip\noindent

For notational simplicity we shall write $\Phi(f)$ for $\Phi(f;\rho,\psi)$. For any $f_1,\,f_2\in C[-1,1]$ we define the convolution integral
$$
(f_1*f_2)(x)=\int_{0}^{x}f_1(x-y)f_2(y)\,dy.\eqno (5.2)
$$

\medskip\noindent
It is direct to verify that
$$
\ps*\varphi_{s'}=\al\frac{\ps-\varphi_{s'}}{s-s'}\qquad (s\neq s').\eqno(5.3)
$$

\medskip\noindent

Write
$$
\delta=R^{-9/10}\qquad\deo=R^{-8/9}.\eqno (5.4)
$$
Fix a constant $d\ge 425$ (in fact, the choice $d=425$ is admissible).  We define two weight functions as follows.
$$
\vpo(x)=\delta^{3-d}(1+\delta-x^2)^{d},\eqno (5.5)
$$
and
$$
\vpt(x)=2d\delta^{3-d}(1+\delta+x^2)(1+\delta-x^2)^{d-2}. \eqno (5.6)
$$

\medskip\noindent
For any $f_1,f_2\in C^1[-1,1]$ we define the inner product
$$
<f_1,\,f_2>=\int_{-1}^{1}\big[f_1'(x)\bar{f}_2'(x)\vpo(x)+f_1(x)\bar{f}_2(x)\vpt(x)\big]\,dx.\eqno (5.7)
$$
We also define the norm
$$
\|f\|=\sqrt{<f,\,f>}\qquad \big(f\in C^1[-1,1]\big).\eqno (5.8)
$$
\medskip\noindent

Since $\vpo(x)\gg\de^3$ for $|x|\le 1$ and $\vpo(x)\gg R^2$ for $|x|\le 1-\deo$, it follows that
$$
\int_{-1}^{1}f'(x)|^2\,dx\ll\de^{-3}\|f\|^2,\qquad (f\in C^1[-1,1]),\eqno (5.9)
$$
and
$$
\int_{-1+\deo}^{1-\deo}f'(x)|^2\,dx\ll R^{-2}\|f\|^2,\qquad (f\in C^1[-1,1]).\eqno (5.10)
$$
The relation (5.10) also holds with $f'(x)$ replaced by $f(x)$. Thus, for any $f\in C^1[-1,1]$ there exists a number $x_1\in[-1+\deo,1-\deo]$ such that $f(x_1)\ll R{-1}\|f\|$, so, by the Cauchy inequality and (5.9),
$$
\max_{|x|\le 1}\{|f(x)|\}\ll R{-1}\|f\|+\int_{-1}^{1}f'(y)|\,dy\ll\de^{-1}\|f\|.\qquad (f\in C^1[-1,1]).\eqno (5.11)
$$
Here we have used the inequality
$$
\int_{-1}^{1}\vpo(y)^{-1}\,dy\ll\de^{-2}.
$$
Analogously we have
$$
\max_{|x|\le 1-\deo}\{|f(x)|\}\ll R{-1}\|f\|.\qquad (f\in C^1[-1,1]).\eqno (5.12)
$$
\medskip\noindent

We shall frequently encounter the problem of bounding $\Phi(f)$ or its partial sums on average when the function $f$ depends on $\psi$ and $\rho$. In such a situation the typical techniques of mean value estimates of Dirichlet polynomials with characters fail to apply directly. Our treatment is to separate the terms with $n=1$ in $\Phi(f)$ and transform the remaining sums into simple sums by partial summation. To this end we introduce
$$
\Phi_1(f)=\Phi_1(f;\rho,\psi)=\sum_{1<n\le Q^{3/2}}\frac{\lp(n)\psi(n)}{n^{\rho}}\tho f(1-\al\log n)\eqno (5.13)
$$
$$
\qquad\qquad-\sum_{1<n\le Q}\frac{\lm(n)\bp(n)}{n^{1-\rho}}\tht f(\al\log n-1), 
$$ 
so
$$
\Phi(f)=\vartheta(Q^{4/3})f(1)-\vartheta(Q^{2/3})f(-1)+\Phi_1(f).\eqno (5.14)
$$

\medskip\noindent
{\bf Remark.} By virtue of (4.16), in practice, the factors $\vartheta(Q^{4/3})$ and $\vartheta(Q^{2/3})$ in (5.14) can be replaced by one with negligible errors.

For $b\ge a\ge 0$ we define
$$
\u_+(a,b;s,\psi)=\sum_{Q^a<n\le Q^b}\frac{\lp(n)\psi(n)}{n^s}\tho,
$$
$$
\u_-(a,b;s,\psi)=\sum_{Q^a<n\le Q^b}\frac{\lm(n)\psi(n)}{n^s}\tht.
$$
Note that the term with $n=1$ is excluded from the above sums in the case $a=0$.
Write
$$
\up\rp=\upp\rp+\upm\rp+\up^*\rp,\eqno (5.15)
$$
where
$$
\begin{aligned}
\up_{\pm}\rp=|\u_{\pm}(0,\deo;\rho,\psi)|^2+|\u_{\pm}(\deo,1;\rho,\psi)|^2+&\int_{0}^{\deo}|\u_{\pm}(0,y;\rho,\psi)|^2\vpo(1-y)^{-1}\,dy\\
\,\\
&+\int_{\deo}^{1}|\u_{\pm}(\deo,y;\rho,\psi)|^2\,dy,
\end{aligned}
$$
and
$$
\up^*\rp=\ep^{-2}|\u_{+}(0,\ep;\rho,\psi)|^2+|\u_{+}(1,3/2;\rho,\psi)|^2+\int_{0}^{1/2}|\u_{+}(1+y,3/2;\rho,\psi)|^2\,dy.
$$

\medskip\noindent
Although the expression for $\up\rp$ is complicated, it can be shown that the average order of $\up\rp$ is $O(1)$. The following assumption will prove very convenient at various places.

\medskip\noindent
$$
\Upsilon\rp<(\log R)^2.\eqno (A2)
$$

\medskip\noindent
{\bf Lemma 5.2} {\it Assume (A2) holds. Then for any $f\in C^1[-1,1]$ we have}
$$
|\Phi_1(f)|\ll\|f\|\lr.
$$

\medskip\noindent
{\it Proof.} We can assume $\|f\|>0$ without loss of generality. Let $\eta$ be a positive number to be specified later. We have
$$
\begin{aligned}
\bigg|\sum_{1<n\le Q^{\deo}}\frac{\lp(n)\psi(n)}{n^{\rho}}\tho f(1-\al\log n)\bigg|
\le|\u_+(0,\delta_1;\rho,\psi)||f(1-\deo)|\\
+\bigg|\sum_{1<n\le Q^{\deo}}\frac{\lp(n)\psi(n)}{n^{\rho}}\tho\int_{\al\log n}^{\deo}f'(1-y)\,dy\bigg|.
\end{aligned}
$$
Interchanging the order of integration and summation we see that the second term on the right-hand side above does not exceed
$$
\int_{0}^{\deo}|\u_+(0,y;\rho,\psi)||f'(1-y)|\,dy.
$$
Combining these results with the trivial inequalities
$$
|\u_+(0,\deo;\rho,\psi)||f(1-\deo)|\le\frac{\eta}{4}|f(1-\deo)|^2+\frac{1}{\eta}|\u_+(0,\deo;\rho,\psi)|^2,
$$
$$
|\u_+(0,y;\rho,\psi)||f'(1-y)|\le\frac{\eta}{4}|f'(1-y)|^2\vpo(1-y)+\frac{1}{\eta}\frac{|\u_+(0,y;\rho,\psi)|^2}{\vpo(1-y)},
$$
we get
$$
\sum_{1<n\le Q^{\deo}}\frac{\lp(n)\psi(n)}{n^{\rho}}\tho f(1-\al\log n)\ll\eta|f(1-\deo)|^2+\frac{1}{\eta}|\u_+(0,\delta_1;\rho,\psi)|^2
$$
$$
+\int_{0}^{\deo}\bigg\{\eta|f'(1-y)|^2\vpo(1-y)+\frac{1}{\eta}\frac{|\u_+(0,y;\rho,\psi)|^2}{\vpo(1-y)}\bigg\}\,dy.
$$
Analogously we have
$$
\sum_{Q^{\deo}<n\le Q}\frac{\lp(n)\psi(n)}{n^{\rho}}\tho f(1-\al\log n)\ll\eta|f(0)|^2+\frac{1}{\eta}|\u_+(\deo,1;\rho,\psi)|^2
$$
$$
+\int_{\deo}^{1}\bigg\{\eta|f'(1-y)|^2\vpo(1-y)+\frac{1}{\eta}\frac{|\u_+(\delta_1,y;\rho,\psi)|^2}{\vpo(1-y)}\bigg\}\,dy,
$$
and
$$
\sum_{Q^<n\le Q^{3/2}}\frac{\lp(n)\psi(n)}{n^{\rho}}\tho f(1-\al\log n)\ll\eta|f(0)|^2+\frac{1}{\eta}|\u_+(1,3/2;\rho,\psi)|^2
$$
$$
+\int_{1}^{3/2}\bigg\{\eta|f'(1-y)|^2\vpo(1-y)+\frac{1}{\eta}\frac{|\u_+(y,3/2;\rho,\psi)|^2}{\vpo(1-y)}\bigg\}\,dy.
$$

\medskip\noindent
The second sum in the expression for $\Phi_1(f)$ is analogously estimated. Note that $\lm(n)$ and 
$\vartheta(x)$ are real, $\overline{1-\rho}=\rho$ by Proposition 3.1, and $\vpo(y-1)=\vpo(1-y)$. Thus we have
$$
\sum_{1<n\le Q^{\deo}}\frac{\lm(n)\bp(n)}{n^{1-\rho}}\tht f(\al\log n-1)\ll
\eta|f(-1+\deo)|^2+\frac{1}{\eta}|\u_-(0,\deo;\rho,\psi)|^2
$$
$$
+\int_{0}^{\deo}\bigg\{\eta|f'(-1+y)|^2\vpo(1-y)+\frac{1}{\eta}\frac{|\u_-(0,y;\rho,\psi)|^2}{\vpo(1-y)}\bigg\}\,dy,
$$
and
$$
\sum_{Q^{\deo}<n\le Q}\frac{\lm(n)\bp(n)}{n^{1-\rho}}\tht f(\al\log n-1)
\ll\eta|f(-1+\deo)|^2+\frac{1}{\eta}|\u_-(\deo,1;\rho,\psi)|^2
$$
$$
+\int_{\deo}^{1}\bigg\{\eta|f'(-1+y)|^2\vpo(1-y)+\frac{1}{\eta}\frac{|\u_-(\delta_1,y;\rho,\psi)|^2}{\vpo(1-y)}\bigg\}\,dy.
$$
Combining the above results with (5.12) and the simple relation
$$
\vpo(1-y)>1 \qquad\text{if}\quad \deo\le y\le 3/2,
$$
we conclude that
$$
\Phi_1(f)\ll\eta\|f\|^2+\frac{1}{\eta}\up\rp.
$$
Hence, choosing $\eta=\|f\|^{-1}\lr$ and applying (A2) we complete the proof. $\Box$
\medskip\noindent

\section *{6. Solutions to Boundary Value Problems}
\medskip\noindent

In this section, we introduce a number of functions which are solutions to certain Sturmian boundary value problems. All the functions in this section are defined on $[-1,1]$.

Let
$$
g_1(x)=(2d)^{-1}\delta^{d-2}(1+\delta-x^2)^{-d}.
$$
It is direct to verify that $g_1$ satisfies the differential equation
$$
[g_1'\vpo]'-g_1\vpt=0\eqno (6.1)
$$
with the boundary condition
$$
g_1'(1)\vpo(1)=-g_1'(-1)\vpo(-1)=1. \eqno (6.2)
$$
Write
$$
\aa=4d^2\delta^{1-d}\int_{0}^{1}(1+\delta-x^2)^d\,dx+2d\delta^2,
$$
and
$$
g_2(x)=2d(\aa\delta)^{-1}(1+\delta-x^2)^{-d}\int_{0}^{x}(1+\delta-y^2)^d\,dy.
$$
It is direct to verify that the function $g_2$ satisfies the same differential equation as $g_1$ and the boundary condition
$$
g_2'(1)\vpo(1)=g_2'(-1)\vpo(-1)=1.\eqno (6.3)
$$
By (6.1), (6.2), (6.3) and partial integration we get for any $f\in C^1[-1,1]$
$$
<f,\,g_1>=f(1)+f(-1),\eqno (6.4)
$$
$$
<f,\,g_2>=f(1)-f(-1).\eqno (6.5)
$$
In particular, we have
$$
\|g_j\|^2=2g_j(1),\qquad j=1,2.
$$
It follows, by routine calculation, that
$$
\|g_j\|^2=d^{-1}\delta^{-2}+O(\delta^{d-1})\qquad j=1,2.\eqno (6.6)
$$

Now we introduce the solution to an inhomogeneous differential equation. Let $\ep$ be a real number satisfying
$$
\ep=R^{-11/12}+O(R^{-1}),\qquad\ep\equiv 0 (\mod 2\pi). \eqno (6.7)
$$
Let $\t_0\in C^2[-1,1]$ such that
$$
0\le\t_0(x)\le 1, \qquad \t'_0(x)\ll\deo^{-1},
$$
$$
\t_0(x)=
\begin{cases} 1,\qquad&\text{if}\quad 1-\deo\le|x|\le 1,\\
0,\qquad&\text{if}\quad |x|\le 1-2\deo.
\end{cases}
$$
Let
$$
\t_1(x)=\frac{\sin\big(R(1+\ep-x)\big)}{1+\ep-x}\t_0(x)\qquad\text{if}\quad x\ge 0,
$$
$$
\t_1(x)=\frac{\sin\big(R(-1+\ep-x)\big)}{-1+\ep-x}\t_0(x)\qquad\text{if}\quad x<0.
$$
Let $\t_2\in C[-1,1]$ be such that 
$$
\t_2(x)=0\qquad\text{unless}\quad |x-(-1+\ep)|<R^{-10},
$$
$$\t_2(x)\ll R^{10},\qquad\int_{-1}^{1}\t_2(x)\,dx=1.
$$
Let 
$$
\t(x)=\t_1(x)-\t_2(x).
$$

Let $g_3(x)$ be the (unique) solution to the boundary value problem
$$
\big[\vpo g'\big]'-\vpt g=-\t,
$$
$$
g'(1)=g'(-1)=0.
$$
It follows, by partial integration, that
$$
<f,g_3>=\int_{-1}^{1}f(y)\t(y)\,dy\qquad (f\in C^1[-1,1]).\eqno (6.8)
$$
By the method of variation of constants we have (see W. Walter [W], for example)
$$
g_3(x)=-\tg(x)+\frac12\big[\tg'(1)-\tg'(-1)\big]\vpo(1)g_1(x)+\frac12\big[\tg'(1)+\tg'(-1)\big]\vpo(1)g_2(x),
$$
where
$$
\tg(x)=\delta^{d-3}(1+\delta-x^2)^{-d}\int_{0}^{x}\bigg\{\int_{y}^{x}(1+\delta-z^2)^{d}\,dz\bigg\}(1+\delta-y^2)^{-d}\t(y)\,dy.
$$
Note that $(\delta\ep R)^{-1}\sim R^{49/60}$. By partial integration and routine estimation we obtain the following

\medskip\noindent
{\bf Lemma 6.1} {\it We have}
$$
\|\tg\|\ll R^{49/60}.
$$

\medskip\noindent

\section *{7. Approximation to $\pia$ }
\medskip\noindent

Recall the definition of $T\rp$ given in Section 4. Let $\ss\rp$ denote the complex vector space spanned by the exponential functions $\pt$ with $\theta\in T\rp$. For notational simplicity we shall write $T$ and $\ss$ for $T\rp$ and $\ss\rp$ respectively.

Note that $\pia(x)\equiv 1$. By the Jensen formula, it can be shown that there exists a function $f\in\ss$ such that
$$
\|\pia-f\|\ll\de.
$$
(see P. Borwein and T. Erd\'{e}lyi [BE, Chapter 4], for example). This result can not be essentially improved. On the other hand, if one allows $f\in\ss+<\ps>$ with $s\ll\al$ and $s\notin T$, then much better approximation to $\pia$ can be obtained. The main result of this section is the following

\medskip\noindent
{\bf Lemma 7.1.} {\it Assume (A2) holds. There exist functions 
$$
h=\sum_{\theta\in T}\aht\pt\in\ss
$$
and $k\in C^1[-1,1]$ such that}
$$
\aht\ll 1,\eqno(7.1)
$$
$$
$$
$$
h(x)\ll(\lr)^2,\qquad h'(x)\ll R\lr\qquad\text{if}\quad|x|\le 2,\eqno (7.2)
$$
$$
$$
$$
B(s):=\sum_{\theta\in T}\frac{\aht}{s-\theta}\ll\de\qquad\text{if}\quad s\in\c,\eqno(7.3)
$$
$$
$$
$$
<k,g_3>=0,\qquad \|k\|\ll\delta,\eqno (7.4)
$$
{\it and such that}
$$
\|r\|=o(\de), \eqno(7.5)
$$
where
$$
r(x)=\pia(x)-h(x)-k(x).\eqno (7.6)
$$

\medskip\noindent
{\it Proof.} We give a sketch. Write $T_j$ for $T_j\rp$, $j=1,2$, and let $\ss_j$ be the complex vector space spanned by the functions $\pt$ with $\theta\in T_j$, $j=1,2$. Let $\r^*$ be a rectangle with vertices at 
$$
\pm\omega+O(\al)\pm i\omega(\lr)^{-1}
$$
such that $\k_+\rsp^{-1}\ll 1$ on the horizontal segments of $\r^*$. Write
$$
\tr_1(x)=\tp\int_{\r^*}\frac{\ps(x)}{Q^s\k_+\rsp}\frac{\exp\big[s^2(\lr)^2/\omega^2\big]}{s}\,ds,
$$
$$
\tr_2(x)=\tp\int_{\r^*}\frac{\ps(x)}{Q^s\k_+\rsp}\frac{\exp\big[s^2(\lr)^2/\omega^2\big]}{s-\al}\,ds.
$$
Both of $\tr_1(x)$ and $\tr_2(x)$ are $\sim 0$ unless $x\sim\pm 1$, while the difference
$$
\tr_1(x)-\tr_2(x)=-\frac{\al}{2\pi i}\int_{\r^*}\frac{\ps(x)}{Q^s\k_+\rsp}\frac{\exp\big[s^2(\lr)^2/\omega^2\big]}{s(s-\al)}\,ds
$$
is much closer to zero. Because of the factor $\exp\big[s^2(\lr)^2/\omega^2\big]$, the integrals on the horizontal segments in the expressions for $\tr_j(x)$ can be neglected. On the other hand, by Lemma 4.3 and (A2), we have $\k_+\rsp\sim 1$ if $s$ lies on the right side of $\r^*$, and $\k_+\rsp\sim -Q^{-2s}$ if $s$ lies on the left side of $\r^*$. It follows that
$$
\tr_j(x)\sim\tp\int_{(\omega)}Q^{(x-1)s}\frac{\exp\big[s^2(\lr)^2/\omega^2\big]}{s}\,ds\qquad(s\sim 1),\eqno (7.7)
$$
and
$$
\tr_j(x)\sim\tp\int_{(-\omega)}Q^{(x+1)s}\frac{\exp\big[s^2(\lr)^2/\omega^2\big]}{s}\,ds\qquad(s\sim -1). \eqno (7.8)
$$
Furthermore, by the residue theorem, there exists a function $h_1\in\ss_1$ such that
$$
\tr_1=\k_+\rp^{-1}\pia+h_1. \eqno (7.9)
$$

\medskip\noindent
By virtue of (7.7) and (7.8), we can construct a function 
$$
h_2=\sum_{\theta\in T_2}\aft\pt\in\ss_2
$$
such that
$$
\aft\ll 1,
$$
$$
$$
$$
\sum_{\theta\in T_2}\frac{\aft}{\theta}=-\k_+\rp, \eqno (7.10)
$$
and such that the function
$$
k(x):=(h_2*\tr_2)(x)
$$
satisfies 
$$
<k,\,g_3>=0\qquad\text{and}\quad\|k\|\ll\de.
$$ 
Let 
$$
r(x)=(\tr_1*h_2)(x)-k(x).
$$ 
It can be shown that 
$$
\|r\|=o(\de).
$$
On the other hand, by (5.3), (7.9) and (7.10), we see that there exists a function
$h\in\ss$ such that
$$
(\tr_1*h_2)(x)=\pia(x)-h(x).
$$
Combining these results we get the assertion. $\Box$

\medskip\noindent
{\bf Remark.} The estimate (7.5) will be sufficient for our purpose; in fact a sharper estimate can be proved.

\section *{8. The Fundamental Inequality: Preliminary }
\medskip\noindent

In this section we assume (A2) holds and start to derive the fundamental inequality described in Section 1. 
\medskip\noindent

Recall that $\ep$ is given by (6.7).
Let $\c_+$ be the right half of the circle $\c$ from $-i\omega$ to $i\omega$, $\c_-$ the left half of $\c$ from $i\omega$ to $-i\omega$. For any $f\in C[-1,1]$ we define three linear functionals as follows.

Let
$$
\begin{aligned}
\Phi^*(f)=\Phi^*(f;\rho,\psi)=&\sum_{Q^{\ep}<n\le Q^{3/2}}\frac{\lp(n)\psi(n)}{n^{\rho}}\tho f(1+\ep-\alpha\log n)\\
\,\\
&-\sum_{1<n\le Q}\frac{\lm(n)\bp(n)}{n^{1-\rho}}\tht f(\alpha\log n-1+\ep).
\end{aligned}
$$
Recall that the convolution integral $*$ is given by (5.2). Let

$$
\begin{aligned}
\Xi(f)=\Xi(f;\rho,\psi)=&\frac{\l^2}{2\pi i}\int_{\c_+}\k_-(\rho+s,\psi)\Phi(f*\ps)Q^{\ep s}\,ds\\
\,\\
&-\frac{e^{-2}\l^2}{2\pi i}\int_{\c_-}\k_-(\rho+s,\psi)\Phi(f*\ps)Q^{(2+\ep) s}\,ds,
\end{aligned}
$$
and let
$$
\Theta(f)=\Theta(f;\rho,\psi)=-(1-e^{-2})f(-1+\ep)+\Xi(f)+\Phi^*(f). 
$$

\medskip\noindent
Write
$$
\hh(x)=h(x+\ep).
$$
 By (7.1) and Lemma 5.1 we get
$$
\Phi(\hh)=\sum_{\theta\in T}Q^{\ep\theta}\aht\Phi(\pt)\ll\l^{-1/3}. \eqno (8.1)
$$

\medskip\noindent
On the other hand, we have
$$
\Phi(\hh)=h(1+\ep)-h(-1+\ep)+\p\rp+\Phi^*(h),\eqno (8.2)
$$
where
\medskip\noindent

$$
\p\rp=\sum_{1<n\le Q^{\ep}}\frac{\lp(n)\psi(n)}{n^{\rho}}h(1+\ep-\al\log n).
$$

\medskip\noindent
We first prove an upper bound for $\p\rp$ under the assumption (A2). By virtue of the proof of Lemma 5.2, we deduce from (7.2) that
$$
\p\rp\ll(\lr)^2\,|\u(0,\ep;\rho,\psi)|+R\lr\int_{0}^{\ep}|\u(0,y;\rho,\psi)|\,dy\ll\deo.
$$
Combining this with (8.1) and (8.2) we get
$$
h(1+\ep)-h(-1+\ep)+\Phi^*(h)\ll\deo.\eqno (8.3)
$$
\medskip\noindent

We claim that
$$
h(1+\ep)-e^{-2}h(-1+\ep)=\Xi(h)+O(\deo).  \eqno (8.4)
$$
Noting that
$$
|1-e^{\pm 2}Q^{-2s}|\gg 1\qquad\text{for}\quad s\in\c_+, 
$$
by Lemma 4.1 and Lemma 4.2 we get
$$
1\ll|\k_{\pm}\rsp|\ll 1\qquad (s\in\c_+).  \eqno (8.5)
$$
On the other hand, noting that 
$$
|1-e^{\pm 2}Q^{-2s}|\gg Q^{-2\sigma}\qquad\text{for}\quad s\in\c_-, 
$$
by Lemma 4.1 and Lemma 4.2 we get
$$
1\ll Q^{-2\sigma}|\k_{\pm}\rsp|\ll 1\qquad (s\in\c_-).  \eqno (8.6)
$$
The proof of (8.4) is divided into five steps.

\bigskip\noindent
{\it Step 1.} We have 
$$
h(1+\ep)-e^{-2}h(-1+\ep)=\sum_{\theta\in T}\aht(1-e^{-2}Q^{-2\theta})Q^{(1+\ep)\theta}.
$$
By (4.13) and (7.1), the terms with $|\theta|>\omega$ in the above sum totally contribute $O(\l^{-1/3})$. (This is a crucial application of the gap assertion (iii) of Proposition 3.1). Hence, by the residue theorem we get
$$
h(1+\ep)-e^{-2}h(-1+\ep)=\tp\int_{\c}B(s)(1-e^{-2}Q^{-2s})Q^{(1+\ep)s}\,ds+O(\l^{-1/3}),\eqno (8.7)
$$
where $B(s)$ is given by (7.3).

\bigskip\noindent
{\it Step 2.}  By Lemma 4.1 we get
$$
1-e^{-2}Q^{-2s}=\frac{\k_-(\rho+s,\psi)\k_+(\rho+s,\psi)}{1-e^{2}Q^{-2s}}+O(\l^{-1/2})\qquad (s\in\c).
$$

\medskip\noindent
We insert this relation into the right-hand side of (8.7) and apply the relation
$$
\frac{1}{1-e^{2}Q^{-2s}}=-\frac{e^{-2}Q^{2s}}{1-e^{-2}Q^{2s}}
$$

\medskip\noindent
getting
$$
h(1+\ep)-e^{-2}h(-1+\ep)=\j_+-\j_-+O(\l^{-1/3}),\eqno (8.8)
$$
where
\medskip\noindent

$$
\j_{+}=\tp\int_{\c_{+}}\frac{\k_-(\rho+s,\psi)\k_+(\rho+s,\psi)}{1-e^{-2}Q^{-2s}}\,B(s)\,Q^{(1+\ep) s}\,ds,
$$
$$
$$
$$
\j_{-}=\frac{e^{-2}}{2\pi i}\int_{\c_{-}}\frac{\k_-(\rho+s,\psi)\k_+(\rho+s,\psi)}{1-e^{2}Q^{2s}}\,B(s)\,Q^{(3+\ep) s}\,ds.
$$

\bigskip\noindent
{\it Step 3.}  We insert the relation
$$
\frac{1}{1-e^2Q^{-2s}}=1+\frac{e^2Q^{-2s}}{1-e^{-2}Q^{-2s}}
$$ 

\medskip\noindent
into the expression for $\j_+$ and then split it into two integrals accordingly.  By (7.3) and routine estimation, we see that the second integral is $O(\deo)$. Hence
$$
\j_+=\tp\int_{\c_+}\k_-(\rho+s,\psi)\k_+(\rho+s,\psi)B(s)Q^{(1+\ep) s}\,ds+O(\deo). 
$$

\medskip\noindent
Analogously, we insert the relation
$$
\frac{1}{1-e^2Q^{2s}}=1+\frac{e^2Q^{2s}}{1-e^2Q^{2s}}
$$ 

\medskip\noindent
into the expression for $\j_-$ and then split it into two integrals accordingly. By (7.3) and routine estimation, we see that the second integral is also $O(\deo)$. Hence
$$
\j_-=\frac{e^{-2}}{2\pi i}\int_{\c_-}\k_-(\rho+s,\psi)\k_+(\rho+s,\psi)B(s)Q^{(3+\ep) s}\,ds+O(\deo). 
$$

\medskip\noindent
Combining these results with (8.8) we get
$$
h(1+\ep)-e^{-2}h(-1+\ep)=\tp\int_{\c_+}\k_-(\rho+s,\psi)\k_+(\rho+s,\psi)B(s)Q^{(1+\ep) s}\,ds \eqno (8.9)
$$
$$
\qquad\qquad -\frac{e^{-2}}{2\pi i}\int_{\c_-}\k_-(\rho+s,\psi)\k_+(\rho+s,\psi)B(s)Q^{(3+\ep) s}\,ds+O(\deo). 
$$

\bigskip\noindent
{\it Step 4.} By Lemma 5.1, we can replace the factor $\k_+\rsp Q^s$ in (8.19) by $\Phi(\ps)$ with an error much smaller than $\deo$. Hence
$$
h(1+\ep)-e^{-2}h(-1+\ep)=\tp\int_{\c_+}\k_-(\rho+s,\psi)\Phi(\ps)B(s)Q^{\ep s}\,ds \eqno (8.10)
$$
$$
\qquad\qquad -\frac{e^{-2}}{2\pi i}\int_{\c_-}\k_-(\rho+s,\psi)\Phi(\ps)B(s)Q^{(2+\ep) s}\,ds+O(\deo). 
$$

\bigskip\noindent
{\it Step 5.} It follows from (5.3) that
$$
h*\ps=-\al\sum_{\theta\in\t}\frac{\aht}{s-\theta}\pt+\al B(s)\ps\qquad (s\in\c),
$$
so
$$
\Phi(h*\ps)=-\al\sum_{\theta\in\t}\frac{\aht}{s-\theta}\Phi(\pt)+\al B(s)\Phi(\ps)\qquad (s\in\c).
$$
The first sum on the right-hand side above is $O(\l^{-1/3})$ by (7.1) and Lemma 5.1. Thus we have
$$
B(s)\Phi(\ps)=\l^2\Phi(h*\ps)+O(\l^{5/3})\qquad (s\in\c). 
$$

\medskip\noindent
Inserting this relation into (8.10) and using routine estimation we get (8.4).
\medskip\noindent

Now, combining(8.4) with (8.3) we get
$$
\Theta(h)\ll\deo.
$$
This, together with (7.6), leads to the following

\medskip\noindent
{\bf Lemma 8.1.} {\it Assume (A2) holds. Then we have}
$$
\Theta(\pia)-\Theta(k)-\Theta(r)\ll\deo.
$$

\section *{9. Upper Bounds for $\Theta(k)$ and $\Theta(r)$}
\medskip\noindent

In this section we assume (A2) holds and derive upper bounds for $\Theta(k)$ and $\Theta(r)$ from the following 

\medskip\noindent
{\bf Lemma 9.1} {\it Assume (A2) holds. For any $f\in C^1[-1,1]$ satisfying $\|f\|\le\de$ we have
$$
\Theta(f)=(1-e^{-2})<f,g_3>+O(\deo)+O(\e_1\rp),
$$
where}
$$
\e_1\rp=\int_{1-2\deo}^{1}\big\{|\x_+(1+\ep-y;\rho,\psi)|^2+|\x_-(1+\ep-y;1-\rho,\bp)|^2\big\}\,dy\eqno (9.1)
$$
$$
+\int_{-1}^{-1+2\deo}\big\{|\x_+(-1+\ep-y;\rho,\psi)|^2+|\x_-(-1+\ep-y;1-\rho,\bp)|^2\big\}\,dy
$$
with the functions $\x_{\pm}$ to be introduced at the end of this section.

\medskip\noindent
{\it Proof.} We begin with estimating $\Phi^*(f)$. By virtue of the proof of Lemma 5.3 with a minor modification we have
$$
\Phi^*(f)\ll\deo. \eqno(9.2)
$$
Here we have used (A2). On the other hand, recalling the definition of $\t_2$ in Section 6, by routine estimation we get
$$
f(-1+\ep)=\int_{-1}^{1}f(x)\t_2(x)\,dx+O(R^{-1}).
$$
\medskip\noindent
It follows that
$$
\Theta(f)=-(1-e^{-2})\int_{-1}^{1}f(x)\t_2(x)\,dx+\Xi(f)+O(\deo).\eqno (9.3)
$$
\medskip\noindent

Now we proceed to estimate $\Xi(f)$. By the Cauchy theorem,
$$
\Xi(f)=\frac{\l^2}{2\pi}\int_{-\omega}^{\omega}(1+e^{-2}Q^{2it})\k_-(\rho+it,\psi)\Phi(f*\varphi_{it})Q^{i\ep t}\,dt.\eqno (9.4)
$$
We have 
$$
\begin{aligned}
\Phi(f*\varphi_{it})=&\sum_{n\le Q^{3/2}}\frac{\lp(n)\psi(n)}{n^{\rho}}\tho (f*\varphi_{it})(1-\al\log n)\\
\,\\
&-\sum_{n\le Q}\frac{\lm(n)\bp(n)}{n^{1-\rho}}\tht (f*\varphi_{it})(\al\log n-1),
\end{aligned}
$$
and
$$
(f*\varphi_{it})(x)=\int_{0}^{x}f(x-y)Q^{iyt}\,dy.\eqno (9.5)
$$
In the case $x<0$ it is convenient to write this function as
$$
(f*\varphi_{it})(x)=-\int_{0}^{-x}f(x+y)Q^{-iyt}\,dy\qquad (x<0).\eqno (9.6)
$$

\medskip\noindent
Hence, interchanging the order of summation and integration we get
$$
\Phi(f*\varphi_{it})=(f*\varphi_{it})(1)-(f*\varphi_{it})(-1)+\v_1(it;\rho,\psi)-\v_2(it;\rho,\psi)+\v_3(it;1-\rho,\bp),
$$
where
$$
\v_1(it;\rho,\psi)=\int_{0}^{1}\w_1(y;\rho,\psi)Q^{iyt}\,dy,\eqno (9.7)
$$
with
$$
\w_1(y;\rho,\psi)=\sum_{1<n\le Q^{1-y}}\frac{\lp(n)\psi(n)}{n^{\rho}}\tho f(1-y-\al\log n),\eqno (9.8)
$$
$$
$$
$$
\v_2(it;\rho,\psi)=\int_{0}^{1/2}\w_2(y;\rho,\psi)Q^{-iyt}\,dy,\eqno (9.9)
$$
with
$$
\w_2(y;\rho,\psi)=\sum_{Q^{1+y}<n\le Q^{3/2}}\frac{\lp(n)\psi(n)}{n^{\rho}}\tho f(1+y-\al\log n),\eqno (9.10)
$$
and 
$$
\v_3(it;1-\rho,\bp)=\int_{0}^{1}\w_3(y;1-\rho,\bp)Q^{-iyt}\,dy,\eqno (9.11)
$$
with
$$
\w_3(y;1-\rho,\bp)=\sum_{1<n\le Q^{1-y}}\frac{\lm(n)\bp(n)}{n^{1-\rho}}\tht f(-1+y+\al\log n).\eqno (9.12)
$$
We insert these relations into (9.4) and note that $\k_-(\rho+it,\psi)\ll 1$ if $|t|\le\omega$ by Lemma 4.2, getting
$$
\Xi(f)=\i(f;\rho,\psi)+O\big(\e^*(f;\rho,\psi)\big)+O(Q^{-2007}),\eqno (9.13)
$$
where
$$
\i(f;\rho,\psi)=\frac{\l^2}{2\pi}\int_{-\omega}^{\omega}Q^{i\ep t}(1+e^{-2}Q^{2it})\k_-(\rho+it,\psi)\big[(f*\varphi_{it})(1)-(f*\varphi_{it})(-1)\big]\,dt\eqno (9.14)
$$
$$
\e^*(f;\rho,\psi)=\l^2\int_{-\omega}^{\omega}\big[|\v_1(it;\rho,\psi)|+|\v_2(it;\rho,\psi)|+|\v_3(it;1-\rho,\bp)|\big]\,dt.\eqno (9.15)
$$
(See the remark below (5.14)).
\medskip\noindent

First we estimate $\e^*(f;\rho,\psi)$. Since
$$
|\v_1(it;\rho,\psi)|\le \frac{1}{4 R^2}+R^2|\v_1(it;\rho,\psi)|^2,
$$
and, by a change of variable and the Plancherel theorem,
$$
\l^2\int_{-\omega}^{\omega}|\v_1(it;\rho,\psi)|^2\,dt\ll\int_{0}^{1}|\w_1(y;\rho,\psi)|^2\,dy,
$$
it follows that
$$
\l^2\int_{-\omega}^{\omega}|\v_1(it;\rho,\psi)|\,dt\ll R^{-1}+R^2\int_{0}^{1}|\w_1(y;\rho,\psi)|^2\,dy.\eqno (9.16)
$$
Analogously we have
$$
\l^2\int_{-\omega}^{\omega}|\v_2(it;\rho,\psi)|\,dt\ll R^{-1}+R^2\int_{0}^{1/2}|\w_2(y;\rho,\psi)|^2\,dy,\eqno (9.17)
$$
and
$$
\l^2\int_{-\omega}^{\omega}|\v_3(it;1-\rho,\bp)|\,dt\ll R^{-1}+R^2\int_{0}^{1}|\w_3(y;1-\rho,\bp)|^2\,dy.\eqno (9.18)
$$
We estimate the integral on the right-hand side of (9.16) as follows. Assume $\deo\le y\le 1$. Since (see the proof of Lemma 5.2)
$$
|w_1(y,\rho,\psi)|\le|\u_+(0,1-y;\rho,\psi)||f(0)|+\int_{0}^{1-y}|\u_+(0,z;\rho,\psi)||f'(1-y-z)|\,dz,
$$
it follows, by the Cauchy inequality and Lemma 5.2, that
$$
|w_1(y,\rho,\psi)|^2\ll R^{-1}|\u_+(0,1-y;\rho,\psi)|^2+R^{-1}\int_{0}^{1-y}|\u_+(0,z;\rho,\psi)|^2\,dz.
$$
Hence
$$
\int_{\deo}^{1}|w_1(y,\rho,\psi)|^2\,dy\ll R^{-2}\int_{0}^{1}|\u_+(0,y;\rho,\psi)|^2\,dy \eqno (9.19)
$$

\medskip\noindent
Now assume $0\le y\le\deo$. Note that
$$
\begin{aligned}
|w_1(y,\rho,\psi)|&\le|\u_+(0,\deo;\rho,\psi)||f(1-\deo-y)|+|\u_+(\deo,1-y;\rho,\psi)||f(0)|\\
\,\\
+&\int_{0}^{\deo}|\u_+(0,z;\rho,\psi)||f'(1-y-z)|\,dz+\int_{\deo}^{1-y}|\u_+(\deo,z;\rho,\psi)||f'(1-y-z)|\,dz,
\end{aligned}
$$
and $\vpo(1-y-z)^{-1}\le\vpo(1-z)^{-1}$ if $0\le z\le\deo$. Hence, by the Cauchy inequality, (5.10) and (5.12) we get
$$
|w_1(y,\rho,\psi)|^2\ll R^{-3}|\u_+(0,\deo;\rho,\psi)|^2+R^{-1}|\u_+(\deo,1-y;\rho,\psi)|^2
$$
$$
+R^{-2}\int_{\deo}^{1-y}|\u_+(\deo,z;\rho,\psi)|^2\,dz
+\de^{2}\int_{0}^{\deo}|\u_+(0,z;\rho,\psi)|^2\vpo(1-z)^{-1}\,dz.
$$
Hence
$$
\int_{0}^{\deo}|w_1(y,\rho,\psi)|^2\,dy\ll R^{-3}|\u_+(0,\deo;\rho,\psi)|^2+R^{-2}\int_{\deo}^{1}|\u_+(\deo,z;\rho,\psi)|^2\,dz\eqno (9.20)
$$
$$+\deo\de^2\int_{0}^{\deo}|\u_+(0,z;\rho,\psi)|^2\vpo(1-z)^{-1}\,dz. 
$$
Combining (9.19) and (9.20) with (9.16) and noting that $\deo\de^2R^2\sim R^{-31/45}$ we get
$$
\l^2\int_{-\omega}^{\omega}|\v_1(it;\rho,\psi)|\,dt\ll R^{-1}+R^{-31/45}\up\rp.
$$
The same estimate holds for the integrals of $|\v_2(it;\rho,\psi)|$ and $|\v_3(it;1-\rho,\bp|$ by (9.17) and (9.18). Hence, By (A2) we get
$$
\e^*(f;\rho,\psi)\ll R^{-2/3}. \eqno (9.21)
$$

Now we proceed to evaluate $\i(f;\rho,\psi)$. Assume $|t|\le\omega$.
Note that $f(y)\ll R^{-3/2}$ if $|y|\le 1-\deo$ by (5.12), and recall the definition of $\t_0(y)$  in Section 6. 
By (9.5) and (9.6) we get
$$
(f*\varphi_{it})(1)-(f*\varphi_{it})(-1)=\int_{-1}^{0}f(y)\t_0(y)Q^{-it(1+y)}\,dy \eqno (9.22)
$$
$$
+\int_{0}^{1}f(y)Q^{it(1-y)}\,dy
+O(R^{-3/2}).
$$
On the other hand, we write the result (4.19) on $\k_-\rsp$ in Lemma 4.3 as
$$
\begin{aligned}
\k_-(\rho+it,\psi)=1-Q^{-2it}&+\sum_{1<n\le Q^{3/2}}\frac{\lm(n)\psi(n)}{n^{\rho+it}}\tho\\
\,\\
&-Q^{-2it}\sum_{1<n\le Q}
\frac{\lp(n)\bp(n)}{n^{1-\rho-it}}\tht+O(\l^{-3/5}),
\end{aligned}
$$
so,
$$
\begin{aligned}
(1+e^{-2}Q^{2it})\k_-(\rho+it,\psi)=&1-e^{-2}+e^{-2}Q^{2it}-Q^{-2it}\\
\,\\
&+(1+e^{-2}Q^{2it})\sum_{1<n\le Q^{3/2}}\frac{\lm(n)\psi(n)}{n^{\rho+it}}\tho\\
\,\\
&-(e^{-2}+Q^{-2it})\sum_{1<n\le Q}
\frac{\lp(n)\bp(n)}{n^{1-\rho-it}}\tht+O(\l^{-3/5}).
\end{aligned}
$$
Combining this with the relation
$$
\l^2\int_{-\omega}^{\omega}Q^{itz}\,dt=\frac{2\sin(Rz)}{z},
$$
we get for $|z|\le 2\deo$
$$
\frac{\l^2}{2\pi}\int_{-\omega}^{\omega}Q^{itz}(1+e^{-2}Q^{2it})\k_-(\rho+it,\psi)\,dt=(1-e^{-2})\frac{\sin(Rz)}{\pi z}\eqno (9.23)
$$
$$
\qquad\qquad\qquad+\x_+(z;\rho,\psi)-\x_-(z;1-\rho,\bp)+O(1),
$$
where
$$
\x_+(z;\rho,\psi)=\sum_{1<n\le Q^{3/2}}\frac{\lm(n)\psi(n)}{n^{\rho}}\tho X_+(z-\al\log n)
$$
with
$$
X_+(x)=\frac{1}{\pi}\frac{\sin(Rx)}{x}+\frac{1}{\pi e^2}\frac{\sin(R(2+x))}{2+x},
$$
and
$$
\x_-(z;1-\rho,\bp)=\sum_{1<n\le Q}\frac{\lp(n)\bp(n)}{n^{1-\rho}}\tht X_-(z+\al\log n)
$$
with
$$
X_-(x)=\frac{1}{\pi e^2}\frac{\sin(Rx)}{x}+\frac{1}{\pi}\frac{\sin(R(-2+x))}{-2+x}.
$$
Recall the definition of $\t_1(y)$ in Section 6, and note that $f(y)\ll 1$ by (5.11). Hence, by (9.23), (9.22) and (9.14) we get
$$
\i(f;\rho,\psi)=(1-e^{-2})\int_{-1}^{1}f(y)\t_1(y)\,dy+O\big(\e_1^*\rp\big)+O(\deo),\eqno (9.24)
$$
where
$$
\begin{aligned}
\e_1^*\rp=&\int_{1-2\deo}^{1}\big\{|\x_+(1+\ep-y;\rho,\psi)|+|\x_-(1+\ep-y;1-\rho,\bp)|\big\}\,dy\\
\,\\
&+\int_{-1}^{-1+2\deo}\big\{|\x_+(-1+\ep-y;\rho,\psi)|+|\x_-(-1+\ep-y;1-\rho,\bp)|\big\}\,dy.
\end{aligned}
$$
By a change of variable and the elementary inequality $A\ll 1+|A|^2$ we get
$$
\e_1^*\rp\ll\deo+\e_1\rp.
$$
Hence, by (9.24), (9.21) and (9.13),
$$
\Xi(f)=(1-e^{-2})\int_{-1}^{1}f(y)\t_1(y)\,dy+O(\deo)+O\big(\e_1\rp\big).
$$
Combining this with (9.3) and (6.8) we complete the proof. $\Box$
\medskip\noindent

By Lemma 9.1 and (7.4) we get
$$
\Theta(k)\ll\deo+\e_1\rp.
$$
On the other hand, noting that $(\ep R)^{-1}\sim R^{-1/12}$, by Lemma 9.1, Lemma 6.1, (7.4), (7.5) and the Cauchy-Schwarz inequality we get
$$
\Theta(r)=o( R^{-1/12})+O\big(\e_1\rp\big).
$$
Combining these two relations with Lemma 8.1 we conclude the following

\medskip\noindent
{\bf Lemma 9.2} {\it Assume (A2) holds. Then we have}
$$
\Theta(\pia)=o( R^{-1/12})+O\big(\e_1\rp\big).
$$

\medskip\noindent

\section *{10. Asymptotic Expression for $\Theta(\pia)$}
\medskip\noindent

The aim of this section is to show, apart from some error terms whose average orders will be estimated in Section 12, that
$$
\Theta(\pia)\sim -(1-e^{-2})U_(\ep),
$$
where
$$
U_(\ep)=\tp\int_{\c_-}Q^{\ep s}\,\frac{ds}{s}.
$$
\medskip\noindent

We begin with calculating $\Xi(\pia)$. Write 
$$
U_{\pm}(x)=\tp\int_{\c_{\pm}}\frac{Q^{xs}}{s}\,ds\qquad\qquad(x\in\mathbf{R}).
$$
By the residue theorem we get
$$
U_+(x)+U_-(x)=1.\eqno(10.1)
$$

\medskip\noindent
By (10.1) and routine estimation we have
$$
U_+(x)=1+O((1+Rx)^{-1}),\qquad U_-(x)\ll (1+Rx)^{-1}\qquad (x\ge 0), \eqno (10.2)
$$
and
$$
U_-(x)=1+O((1+R|x|)^{-1}),\qquad U_+(x)\ll (1+R|x|)^{-1}\qquad ( x\le 0). \eqno (10.3)
$$

\medskip\noindent
 By (5.3) we get
$$
\Phi(\pia*\ps)=\frac{\al}{s}\big[-\Phi(\pia)+\Phi(\ps)\big].
$$
Hence
$$
\Xi(\pia*\ps)=-\Phi(\pia)\big[\Xi_1-e^{-2}\Xi_2\big]+\Xi_3-e^{-2}\Xi_4,\eqno (10.4)
$$
where
$$
\Xi_1=\tp\int_{\c_+}\k_-\rsp\, Q^{\ep s}\,\frac{ds}{s},\eqno (10.5)
$$
$$
\Xi_2=\tp\int_{\c_-}\k_-\rsp\, Q^{(2+\ep) s}\,\frac{ds}{s},\eqno (10.6)
$$
$$
\Xi_3=\tp\int_{\c_+}\k_-\rsp \,\Phi(\ps)\,Q^{\ep s}\,\frac{ds}{s}\eqno (10.7)
$$
and
$$
\Xi_4=\tp\int_{\c_-}\k_-\rsp\,\Phi(\ps)\,Q^{(2+\ep) s}\,\frac{ds}{s}.\eqno (10.8)
$$

\medskip\noindent
To evaluate $\Xi_1$ we insert the formula (4.19) into (10.7), and estimate the contribution from the error term to the integral trivially.. Hence, integrating term by term and separating the terms with $n=1$ in the first sum we get
$$
\Xi_1=U_+(\ep)+\y_1(\rho,\psi)-\y_2(1-\rho,\bp)+O(\l^{-1/2}) \eqno (10.9)
$$
where
$$
\y_1\rp=\sum_{1<n\le Q^{3/2}}\frac{\lm(n)\psi(n)}{n^{\rho}}\tho U_+(\ep-\al\log n), \eqno (10.10)
$$
$$
\y_2\rbp=\sum_{n\le Q}
\frac{\lp(n)\bp(n)}{n^{1-\rho}}\tht U_+(-2+\ep+\al\log n). \eqno (10.11)
$$

\medskip\noindent
Analogously we have
$$
\Xi_2=-U_-(\ep)+\y_3(\rho,\psi)-\y_4(1-\rho,\bp)+O \eqno (10.12)
$$
where
$$
\y_3\rp=\sum_{n\le Q^{3/2}}\frac{\lm(n)\psi(n)}{n^{\rho}}\tho U_-(2+\ep-\al\log n), \eqno (10.13)
$$
$$
\y_4\rbp=\sum_{1<n\le Q}\frac{\lp(n)\bp(n)}{n^{1-\rho}}\tht U_-(\ep+\al\log n). \eqno (10.14)
$$
Combining these results we get
$$
\Xi_1-e^{-2}\Xi_2=1-(1-e^{-2})U_-(\ep)+O(\ep)+O(\e_2\rp), \eqno (10.15)
$$
where
$$
\e_2\rp=\ep^{-1}\big[|\y_1\rp|^2+|\y_2\rbp|^2+|\y_3\rp|^2+|\y_4\rbp|^2\big]. \eqno (10.16)
$$
\medskip\noindent

The evaluations of $\Xi_3$ and $\Xi_4$ are simpler. By Lemma 5.1, Lemma 4.2 and Lemma 4.1 we get for $s\in\c$
$$
\k_-\rsp\Phi(\ps)=Q^s-\big(e^{-2}+e^{2}\big)Q^{-s}+Q^{-3s}+O(\l^{-3/10}).
$$
Inserting this relation into (10.7) and (10.8) and applying (10.2) and (10.3) we get
$$
\Xi_3=1+O(R^{-1}), \eqno (10.17)
$$
and
$$
\Xi_4=1+O(R^{-1}). \eqno (10.18)
$$

\medskip\noindent
Combining (10.15), (10.17) and (10.18) with (10.4) we conclude that
$$
\Xi(\pia)=-\Phi(\pia)+(1-e^{-2})U_-(\ep)\Phi(\pia)+1-e^{-2}
\eqno (10.19)
$$
$$
+O(\ep)+O(\e_2\rp).
$$
Note that $\pia(x)\equiv 1$,
and recall the definition of $\Theta$ in Section 8. It follows from (10.19) that
$$
\Theta(\pia)=-\Phi(\pia)+(1-e^{-2})U_-(\ep)\Phi(\pia)+\Phi^*(\pia)\eqno (10.20)
$$
$$
+O(\ep)+O(\e_2\rp).
$$
On the other hand, we have
$$
\Phi^*(\pia)=\Phi(\pia)-\sum_{1<n\le Q^{\ep}}\frac{\lp(n)\psi(n)}{n^{\rho}}\tho.\eqno (10.21).
$$
The sum on the right-hand side above is $O(\ep\lr)$ if (A2) holds. Hence, from (10.20) and (10.21) we derive the following

\medskip\noindent
{\bf Lemma 10.1} {\it Assume (A2) holds. Then we have}
$$
\Theta(\pia)=(1-e^{-2})U_-(\ep)\Phi(\pia)+O(\ep\lr)+O(\e_2\rp).
$$

\section *{11. The Fundamental Inequality: Completion}
\medskip\noindent

It follows from Lemma 5.1, Lemma 4.1 and Lemma 4.2 that
$$
|\Phi(\pia)|\gg 1. \eqno (11.1)
$$
On the other hand, by the Cauchy theorem we get
$$
U_(\ep)=\tp\int_{-i\infty}^{-i\omega}\frac{Q^{\ep s}}{s}\,ds+\tp\int_{i\omega}^{-i\infty}\frac{Q^{\ep s}}{s}\,ds.
$$
Hence, by a change of variable,
$$
U_(\ep)=\frac{1}{\pi}\int_{R}^{\infty}\frac{\sin(\ep y)}{y}\,dy.
$$
Recall that $\ep R\equiv 0 (\mod 2\pi)$. Hence, by partial integration,
$$
|U_(\ep)|\gg (\ep R)^{-1}\sim R^{-1/12}.\eqno (11.2)
$$

\medskip\noindent
Combining (11.1), (11.2) and the results in Section 9 and Section 10 we arrive at

\medskip\noindent
{\bf Proposition 11.1} (The Fundamental Inequality) {\it Assume that $\psi\in\Psi^*$ and $\rho\in\s$. Then we have
$$
\e\rp\gg R^{-1/12},\eqno (11.3)
$$
where}
$$
\e\rp=\frac{1}{R^{1/12}\lr}\up\rp+\e_1\rp+\e_2\rp.
$$

\medskip\noindent
{\it Proof}. If (A2) fails to hold, then we have
$$
\frac{1}{R^{1/12}\lr}\up\rp\ge\frac{\lr}{R^{1/12}}
$$
so (11.3) holds. If (A2) holds, then we have, by (11.1), (11.2), Lemma 9.2 and Lemma 10.1,
$$
\e_1\rp+\e_2\rp\gg R^{-1/12}
$$
so (11.3) holds. $\Box$
\medskip\noindent

Let
$$
\e=\sum_{\psi\in\Psi^*}\,\sum_{\rho\in\s}\e\rp,
$$
From Propositon 11.1 we derive the following 

\medskip\noindent
{\bf Corollary 11.2.} {\it We have}
$$
\e\gg R^{-1/12}\sum_{\psi\in\Psi^*}\,\sum_{\rho\in\s}1.
$$
\medskip\noindent

\section *{12. Estimation of $\e$}
\medskip\noindent

In this section we estimate $\e$. Our result rests only on the fact that $\Psi^*$ is a subset of $\Psi$.

We begin with the following two lemmas.

\medskip\noindent
{\bf Lemma 12.1} {\it For any complex numbers $a_n$ we have}
$$
\sum_{\psi\in\Psi^*}\,\sum_{\rho\in\s}\bigg|\sum_{n\le Q^2}\frac{a_n\psi(n)}{n^{\rho}}\bigg|^2\ll
\l^2Q^2\sum_{n\le Q^2}\frac{|a_n|^2}{n}.
$$

\medskip\noindent
{\it Proof.}  Write
$$
P\sp=\sum_{n\le Q^2}\frac{a_n\psi(n)}{n^{s}}.
$$
Assume that $\psi\in\Psi^*$. For $\rho\in\s$ we have 
$$
|P\rp|^2\ll\l^2\int_{-\al}^{\al}|P(\rho+it,\psi)|^2\,dt+\al\int_{-\al}^{\al}|P'(\rho+it,\psi)|^2\,dt. \eqno(12.1)
$$
By Proposition 3.1, the line segments $[\rho-i\al,\,\rho+i\al]$ and $[\rho'-i\al,\,\rho'+i\al]$ are disjoint if $\rho,\rho'\in\s$ and $\rho\neq\rho'$,
$$
\bigcup_{\rho\in\s}[\rho-i\al,\,\rho+i\al]\subset[s_0-2i,\,s_0+2i].
$$
Hence, by (12.1),
$$
\sum_{\rho\in\s}|P\rp|^2\ll\int_{D-2}^{D+2}\bigg\{\l^2\bigg|P\bigg(\frac12+it,\psi\bigg)\bigg|^2+\al\bigg|P'\bigg(\frac12+it,\psi\bigg)\bigg|^2\bigg\}\,dt\eqno (12.2)
$$
Note that
$$
P'\sp=-\sum_{n\le Q^2}\frac{a_n(\log n)\psi(n)}{n^{s}}
$$
and $\log n\le 2\l^2$ for $n\le Q^2$. Hence,
by (12.2) and Lemma 2.2  we complete the proof. $\Box$
\medskip\noindent

Recall that $\lambda_{\pm}(n)$ are given by (4.18). From the Euler product representation it follows for $\l\ge 1$ that
$$
|\lp(p^l)|\le\lm(p^l)=p^{l\al}(1-p^{-2\al}).\eqno (12.3)
$$

\medskip\noindent
{\bf Lemma 12.2.} {\it For $0<\eta\le 2$ we have}
$$
\sum_{1<n\le Q^{\eta}}\frac{\lpm(n)^2}{n}\ll\eta^2. \eqno (12.4)
$$

\medskip\noindent
{\it Proof.} The left-hand side of (12.4) does not exceed
$$
\prod_{p\le Q^{\eta}}\bigg\{\sum_{0}^{\infty}\frac{\lpm(p^l)^2}{p^l}\bigg\}-1.
$$
Assume $p\le Q^2$ and $l\ge 1$. By (12.3) we get
$$
\lpm(p^l)\ll\al p^{\al l}\log p,
$$
so
$$
\sum_{1}^{\infty}\frac{\lpm(p^l)^2}{p^l}\ll\frac{(\al\log p)^2}{p}.
$$
Hence
$$
\sum_{p\le Q^{\eta}}\log\bigg\{\sum_{0}^{\infty}\frac{\lpm(p^l)^2}{p^l}\bigg\}\ll\al^2\sum_{p\le Q^{\eta}}\frac{(\log p)^2}{p}\ll\eta^2.
$$
This completes the proof. $\Box$
\medskip\noindent

In the estimation of $\e$ we shall repeatedly apply the simple bound $\vartheta(x)\ll 1$ by (4.15). Recall that $\up\rp$ are given by (5.15) and note that 
$$
\int_{0}^{\deo}\frac{y^2\,dy}{\vpo(1-y)}\ll 1.
$$
By Lemma 12.1 and Lemma 12.2 we get
$$
\sum_{\psi\in\Psi^*}\,\sum_{\rho\in\s}\up\rp\ll\l^2Q^2.\eqno (12.5)
$$

Recall that $\e_1\rp$ is given by (9.1). To estimate the sum of $\e_1\rp$ in $\e$ we need the following consequence of Lemma 12.2
$$
\sum_{1<n<Q^2}\frac{\lambda_{\pm}(n)^2}{n(1+R\al\log n)^2}\ll R^{-2}\log R,\eqno(12.6)
$$ 
which follows by partial summation. Note that
$$
X_{\pm}(z\mp\al\log n)\ll\min\{R\,|z\mp\log n|\}
$$
if $|z|\le 3\deo$ and $n\le Q^{3/2}$. Hence, by Lemma 12.1 and (12.6) we get
$$
\sum_{\psi\in\Psi^*}\,\sum_{\rho\in\s}\e_1\rp\ll\l^2Q^2R^{-1/2}.\eqno (12.7)
$$
\medskip\noindent

Recall that $\e_2\rp$ is given by (10.16). By (10.2) and (10.3) we have $U_+(\ep-\al\log n)\ll 1$ and
$$
U_+(\ep-\al\log n)\ll(1+R\al\log n)^{-1}\qquad\text{if}\quad n>Q^{2\ep}.
$$
Hence, by Lemma 12.1 and (12.6)  we get
$$
\sum_{\psi\in\Psi^*}\,\sum_{\rho\in\s}\e_2\rp\ll\l^2Q^2R^{-1/2}.\eqno (12.8)
$$
\medskip\noindent

Finally, combining (12.5), (12.7) and (12.8) we conclude the following

\medskip\noindent
{\bf Proposition 12.3} {\it We have}
$$
\e\ll\l^2Q^2R^{-1/12}(\log R)^{-1}.
$$
\medskip\noindent

\section *{13. Proof of Theorem 1}
\medskip\noindent

Assume $q\in\q$, so $q$ is square-free and every prime divisor of $q$ is greater than three. Let
$$
q=p_1p_2....p_k
$$
be the prime decomposition of $q$. Then every character $\psi\mod q$ can be uniquely written as
$$
\psi=\psi_1\psi_2....\psi_k,
$$
where $\psi_j$ is a character $\mod p_j$. By the definition of the set $\Psi_q$, $\psi$ is in $\Psi_q$ if and only if $\psi_j$ is primitive for $(p_j,D)=1$, and $\psi_j$ is non-real for $p_j|D$. Thus we have
$$
\sum_{\psi\in\Psi_q}1=\prod_{\substack{p|q\\(p,D)=1}}(p-2)\prod_{p|(q,D)}(p-3).
$$
It follows that
$$
N:=\sum_{\psi\in\Psi}1\gg\sum_{q\in\q}q\bigg(\frac{\varphi(q)}{q}\bigg)^3\gg Q^2. \eqno(13.1)
$$
On the other hand, we have
$$
\sum_{\rho\in\s}1\gg\l^2\qquad(\psi\in\Psi^*).\eqno (13.2)
$$
By (13.1), (13.2) and Proposition 2.7 we conclude the following

\medskip\noindent
{\bf Proposition 13.1.} {\it Assume (A1) holds. Then we have}
$$
\sum_{\psi\in\Psi^*}\,\sum_{\rho\in\s}1\gg\l^2Q^2.
$$
\medskip\noindent

Now, under the assumption (A1), a contradiction is immediately derived from Corollary 11.2, Proposition 12.3 and Proposition 13.1.

\medskip\noindent

The proof of Theorem 1 is complete. 

\medskip\noindent
{\bf Remark.} Under the slightly weaker assumption
$$
L(1,\chi)<c(\log D)^{-17}
$$
with the constant $c$ sufficiently small, it can be shown that
$$
\sum_{\psi\in\Psi^*}1\gg Q^2,
$$
and thus a contradiction can be derived.

\bigskip\noindent
{\bf Acknowledgements}. The author thanks Professor Kenneth Appel, Professor Liming Ge and Professor Eric Grinberg of University of New Hampshire for their continuous encouragements and technical helps.

\begin{flushleft}
\medskip\noindent
\begin{tabbing}
XXXXXXXXXXXXXXXXXXXXXXXXXX*\=\kill
Yitang Zhang\\
Department of Mathematics and Statistics\\
University of New Hampshire\\
Durham, NH 03824\\

E-mail: yitangz@cisunix.unh.edu

\end{tabbing}

\end{flushleft}
\end {document}